\documentclass[11pt,fleqn,twoside]{article}

\usepackage{amsthm}
\usepackage{graphicx}
\usepackage{amssymb}
\makeatletter
\newcommand{\prava}[1]{\small\it
\begin{flushleft}
Copyright \copyright \ 2000 by  #1
\end{flushleft}}

\newcommand{\name}[1]{\begin{flushleft}
                       \LARGE \bf #1
                       \end{flushleft}\vspace{-3mm}}

\newcommand{\Author}[1]{\begin{flushleft}
                       \it #1 \end{flushleft}}

\newcommand{\Adress}[1]{\begin{flushleft}
                       \it #1 \end{flushleft}}

\newcommand{\Date}[1]{\begin{flushleft}
                      \small  \it #1 \end{flushleft}}

\newcommand{\ehkol}{Author \ name}
\newcommand{\ohkol}{Article \ name}

\renewcommand{\@evenhead}{
\hspace*{-3pt}\raisebox{-15pt}[\headheight][0pt]{\vbox{\hbox to \textwidth 
{\thepage \hfil \ehkol}\vskip4pt \hrule}}}
\renewcommand{\@oddhead}{
\hspace*{-3pt}\raisebox{-15pt}[\headheight][0pt]{\vbox{\hbox to \textwidth 
{\ohkol \hfil \thepage}\vskip4pt\hrule}}}
\renewcommand{\@evenfoot}{}
\renewcommand{\@oddfoot}{}

\long\def\@makecaption#1#2{%
  \vskip\abovecaptionskip
  \sbox\@tempboxa{\small \textbf{#1.}\ \ #2}%
  \ifdim \wd\@tempboxa >\hsize
    {\small \textbf{#1.}\ \ #2}\par
  \else
    \global \@minipagefalse
    \hb@xt@\hsize{\hfil\box\@tempboxa\hfil}%
  \fi
  \vskip\belowcaptionskip}

\def\numberwithin#1#2{\@ifundefined{c@#1}{\@nocounterr{#1}}{%
  \@ifundefined{c@#2}{\@nocnterr{#2}}{%
  \@addtoreset{#1}{#2}%
  \toks@\@xp\@xp\@xp{\csname the#1\endcsname}%
  \@xp\xdef\csname the#1\endcsname
    {\@xp\@nx\csname the#2\endcsname
     .\the\toks@}}}}

\renewenvironment{proof}[1][\proofname]{\par
  \normalfont
  \topsep6\p@\@plus6\p@ \trivlist
  \item[\hskip\labelsep\textbf{%
    #1}\@addpunct{.}]\ignorespaces
}{%
  \qed\endtrivlist
}

\newcommand{\resetfootnoterule} {
  \renewcommand\footnoterule{%
  \kern-3\p@
  \hrule\@width.4\columnwidth
  \kern2.6\p@}
}

\makeatother

\numberwithin{equation}{section}
\renewcommand{\qedsymbol}{$\blacksquare$}

\newtheorem{lem}{Lemma}[section]
\newtheorem{thm}{Theorem}[section]
\newtheorem{cor}{Corollary}[section]
\theoremstyle{definition}
\newtheorem{rem}{Remark}
\newtheorem{definition}{Definition}[section]

\newcommand{\df}{\displaystyle}
\newcommand{\eps}{\epsilon}
\newcommand{\om}{\Omega}
\newcommand{\p}{\partial}

\begin{document}

\thispagestyle{empty}
\renewcommand{\ehkol}{B.\ Birnir and N.\ Svanstedt}
\renewcommand{\ohkol}{Existence and Homogenization of the 
Rayleigh-B\'enard Problem}

\begin{flushleft}
\footnotesize \sf
Journal of Nonlinear Mathematical Physics \qquad 2000, V.7, N~2,
\pageref{bir_svan_fp}--\pageref{bir_svan_lp}.
\hfill {\sc Article}
\end{flushleft}

\vspace{-5mm}

\renewcommand{\footnoterule}{}
{\renewcommand{\thefootnote}{}
 \footnotetext{\prava{B.\ Birnir and N.\ Svanstedt}}}

\name{Existence and Homogenization of the Rayleigh-B\'enard Problem}
\label{bir_svan_fp}

\Author{Bj\"{o}rn BIRNIR~$^\dag$ 
        and
        Nils SVANSTEDT~$^{\ddag}$}

\Adress{$^\dag$ Department of Mathematics, 
        University of California Santa 
        Barbara, CA 93106, USA \\
        ~~Email: birnir@math.ucsb.edu, 
        URL: www.math.ucsb.edu/\~{}birnir\\[2mm] 
        ~~The University of Iceland, Science Institute, 
        Dunhaga, Reykjav\'{\i}k 107, Iceland\\[2mm]
        $^\ddag$ Department of Mathematics, 
        University of California Santa 
        Barbara, CA 93106, USA \\[2mm]
        ~~Department of Mathematics, Chalmers University 
        of Technology and\\ 
        ~~G\"oteborg University, S-412 96 G\"oteborg, Sweden\\
        ~~Email: nilss@math.chalmers.se, 
        URL: www.math.chalmers.se/\~{}nilss}

\Date{Received June 23, 1999; Revised December 11, 1999; 
Accepted December 17, 1999}

\begin{abstract}
\noindent
The Navier-Stokes equation driven by heat conduction is studied.
As a prototype we consider Rayleigh-B\'enard convection, in the Boussinesq
approximation. Under a large aspect ratio assumption, which is the case
in Rayleigh-B\'enard experiments with Prandtl number close to one, we prove the
existence of a global strong solution to the 3D Navier-Stokes equation 
coupled with a heat equation, and the existence of a maximal B-attractor.
A  rigorous two-scale limit is obtained by homogenization theory.
The mean velocity field is obtained by averaging
the two-scale limit over the unit torus in the local variable.
\end{abstract}

\section{Introduction}

In this paper we study
the Navier-Stokes equation driven by heat conduction.
Under a large aspect ratio 
assumption (the spatial domain being a thin layer) we prove the existence of a global
strong solution. The aspect ratio $\Gamma$ is the width/height ratio of the experimental
apparatus and $\Gamma \ge 40$
is considered large, see Hu et al. \cite{hu1, hu2}. For existence results for the Navier-Stokes equation in thin domains
we particularly refer to Raugel \cite{rau} and the references therein.
Our contribution is the addition of a heat equation where the heat conduction 
is driving the fluid.
The prototype we have in mind is the classical B\'enard problem and this work
is motivated
by the instabilities of roll-patterns observed in experiments and
simulations. Rayleigh-B\'enard convection is a model for pattern formation
and has been extensively studied,  Busse and Clever \cite{bu2,bu1}
established the stability of straight parallel convection rolls and used
them to explain many experimental observations, for large Prandtl numbers
$P  = {\nu \over \kappa}$ ($\nu$ is the kinematic viscosity and $\kappa$
is the thermal diffusitivity). For low Prandtl numbers, the situation is much 
more complicated and it has long been recognized 
\cite{c, cn, dp, nps2} that mean flows are crucial
in understanding the complex pattern dynamics observed in experiments.
In this paper we establish that this complexity is not caused by singularity
formation but
our ultimate goal is a theoretical understanding and a quantitative capture of
the mean flow. It is known that wave-number distortion, roll curvature and 
the mean flow make straight convection rolls become
unstable \cite{slb,hu3,mave,p}. These 
effects have been
successfully modelled by  Decker and Pesch \cite{dp} and 
simulated. However, the resulting equations contain non-local terms due
to the mean-flow and are theoretically intractable. 
It is our hope that our results will put the study of the contribution of the 
mean flow on a rigorous
mathematical footing,  both simplifying and reaching a better
theoretical understanding in the process. 
The experimental difficulties of measuring weak global flow in the presence of
dominant local roll circulations are formidable, \cite{hu2}, and a theoretical
insight may be crucial in the case of dynamical patterns.

As a preparation, we prove a new statement of 
the classical, see Leray \cite{lr1}, existence of 
a global strong solution
and a global attractor for the three-dimensional Navier-Stokes equation under 
smallness assumptions on the data, see Ladyshenskaya \cite{la1}.
For different results with large forcing compare Foias and Temam \cite{FT87}
and Sell \cite{Se96}.
This new statement and new proof of the theorem are crucial in the
statement and the proof of
the existence of a global
solutions of the Rayleigh-B\'enard problem with a large aspect ratio. The 
existence is proven after an initial time-interval, corresponding to a 
settling-down period in experiments. In experiments in gases ($CO_2$) this
settling-down time is a few hours for experiments that take a few days,
\cite{hu2}. 

We will also prove that the Rayleigh-B\'enard problem has a global attractor.
Then theorems of Milnor \cite{ml}, Birnir and Grauer \cite{bg1}, and Birnir 
\cite{bb2} are used to prove that the Rayleigh-B\'enard problem has a unique 
maximal B-attractor. This attractor has the property that every point attracts
a set (of functions) that is not shy \cite{shy}, or of positive 
infinite-dimensional measure.
The discovery of a spiral-defect chaotic attractor \cite{hu4,hu3,mo1}, 
in a parameter region were previously only straight rolls were known to be 
stable \cite{ch, ec1}, was one of the more startling results in recent
pattern formation theory. The experimentally observed attractors are 
B-attractors, and the spiral-defect chaotic attractor and the straight roll
attractor may be two minimal B-attractors of the unique maximal B-attractor
whose existence we prove.

\subsection{The Boussinesq equations}

We start with the Boussinesq equations
which are two coupled equations for the fluid velocity $u$, the pressure $p$ 
and the 
temperature $T$, 
\begin{equation}
\label{eq:bouss}
\left\{ \begin{array}{l}
\df {{\p{u}\over\p{t}}+(u\cdot{\nabla})u-\nu\Delta{u}+\nabla p =
g\alpha(T-T_2)},\\[2ex]
\df {{\p{T}\over\p{t}}+(u\cdot{\nabla})T-\kappa\Delta{T} = 0,}\\[2ex]
{\rm div}\,u = 0, 
\end{array} \right.\;\;x\in\om,\;t\in{\bf R}^+.
\end{equation}
Here $g$ is the gravitational acceleration and
$\alpha$ is the volume expansion coefficient of the fluid.
Moreover $\nu$ and $\kappa$ are the viscosity and conductivity
coefficients which determine the dimensionless Rayleigh 
$R={{\alpha g (T_1-T_2)h^3}\over{\nu\kappa}}$ and Prandtl
numbers. We assume that $\om$ is a rectangular box
and fix the temperatures at the bottom $T_1$ and top $T_2$ of the box.
The box is heated from below so $T_1>T_2$. We can impose periodic boundary
conditions
for $u$ on the lateral sides of the box. However, as the temperature $T_1$
increases these have to be relaxed due to the presence of a boundary layer. 
The velocity is assumed to vanish on the horizontal surfaces of the box. 
Finally, we must supply the appropriate initial data.

\subsection{The homogenization}

The mathematical framework for the homogenization starts by the 
introduction of a (small) parameter
$\eps>0$
and a scaling of the Navier-Stokes system above 
\begin{equation}
\label{eq:h1}
\left\{ \begin{array}{l}
\df {{\p{u_\eps}\over\p{t}}+(u_\eps\cdot{\nabla})u_\eps
-\eps^{3/2} \nu\Delta{u_\eps}+\nabla p_\eps = g\alpha(T_2-T_\eps)},\\[2ex]
\df {{\p{T_\eps}\over\p{t}}+(u_\eps\cdot{\nabla})T_\eps
-\eps^{3/2} \kappa\Delta{T_\eps} = 0,}\\[2ex]
{\rm div}\,u_\eps = 0, 
\end{array} \right.\;\;x\in\om,\;t\in{\bf R}^+.
\end{equation}
The bulk of the paper will be devoted to proving that there exist unique
functions $u_0,\;p_0$ and $T_0$ such that
\[
\eps^{-1/2}u_\eps\to u_0,\;\;p_\eps\to p_0,\;\;
T_\eps\to T_0,
\]
as $\eps\to 0$
in the appropriate Sobolev spaces.
At a first glance, sending $\eps$ to zero
seems to give the Euler equation, but this is wrong. The limit obtained is
viscid
and the above equation only makes sense for $\eps>0$. However, for any fixed
$\eps>0$ the solutions of the scaled equations have global existence
in two dimensions and the equations possess a smooth global attractor 
in dimensions two or three, see Foias et al. \cite{fmt}, Ladyshenskaya \cite{la1,la2}
and Sell \cite{Se96}.
The equations for the leading order coefficients,
in a power series in $\eps$, are still evaluated at $\eps=0$ and turn out to be
the
Navier-Stokes system
\begin{equation}
\left\{ \begin{array}{l}
\df {{\p{u_0}\over\p{\tau}}+(u_0\cdot{\nabla_y})u_0
-\nu\Delta_y{u_0}+\nabla_y p_1=
g\alpha(T_2-T_0)}-\nabla_x p_0,\\[2ex]
\df {{\p{T_0}\over\p{\tau}}+(u_0\cdot{\nabla_y})T_0
-\kappa\Delta_y{T_0} = 0},\\[2ex]
{\rm div_y}\,u_0 = 0,\;\;{\df {\rm div_x}(\int_{T^n} u_0dy)} = 0, 
\end{array} \right.
\end{equation}
where $x\in\om$, $\;y\in T^n$ and $\tau\in{\bf R}^+$. Here
$y=x/\eps$ is the local spatial variable
and $\tau=t/\sqrt{\eps}$ is the scaled fast time variable.
$T^n$, the unit torus in $y$, is what is referred to as the unit cell in the
terminology of homogenization.
The initial data is so
highly oscillatory that we can assume that
the boundary conditions are periodic in the local variable $y$.
 The Navier-Stokes system (1.3) differs from the
original system (1.2) in that it has an additional forcing term $-\nabla_x p_0$.
This is the obvious influence of the global pressure on the local flow. The
solutions
of the system (1.3) enjoy global (in $\tau$) existence in two dimensions
and the equations possess a smooth (in $y$) global attractor in dimensions
two or three. This, possibly high-dimensional, attractor of the local flow
is the physically relevant quantity for the local flow, except in the strong
turbulence
limit, when long transients may play a role. 

We will show that the solution of (1.3) is the unique two-scales limit of a
sequence
of solutions to the scaled system (1.2). This uses the weak sequential
compactness property
of reflexive Banach spaces and says that any bounded sequence $\{u_\eps\}$ in
say $L^2(\om)$ contains a subsequence, still denoted by $\{u_\eps\}$, such that
for
smooth test functions $\varphi(x,y)$, periodic in $y$
\[
\int_\om
u_\eps(x)\varphi(x,{x\over\eps})dx\to\int_\om\int_{T^n}u_0(x,y)\varphi(x,y)dydx.
\]
The main result (Theorem 7.1) is that if $\{u_\eps\}$ is
a sequence of solutions to the Navier-Stokes system (1.2),
then the so called two-scales limit
$u_0$ is the solution to the local Navier-Stokes system (1.3). Our proof is
based upon
a compactness result which was first proved by Nguetseng \cite{ng} and then
further developed by Allaire \cite{a1,a3,a2}. Moreover, if $u_0$ is a globally defined
unique
solution of the system (1.3), which is the case if $u_0$ lies on the attractor
of (1.3),
then, by uniqueness, the whole sequence $\{u_\eps\}$ two-scale converges to
$u_0$.

The mean field turns out to be
\[
{\overline u}_0(x,{t\over\sqrt{\eps}})=\int_0^{t/\sqrt{\eps}}
\left(\pi (e_n{\overline \theta}_0)\right)(x,s)ds,
\]
where 
$e_n$, $n=2,\,3$, is the unit vector in the vertical direction and
$\pi (e_n{\overline \theta}_0)$ denotes the projection onto the divergence free
part of
$e_n{\overline \theta}_0$,
\[
\pi (e_n{\overline \theta}_0)=-\nabla\times(\Delta^{-1}(\nabla\times
e_n{\overline \theta}_0)).
\]
The mean field is derived from the local Navier-Stokes system (1.3).
It gives the contribution of the
conduction to the small scale flow. We have averaged (in $y$), denoted by
overbar, over the unit cell $T^n$, $n=2,\,3$.
Once we have the local problem (4.3)
the boundary conditions on the local cell can also be relaxed to capture the
contribution
of (global) convection to the mean field. The mean field with the influence of
the
convection taken into account, turns out not surprisingly to satisfy a 
forced Euler's equation
\[
{\p{\overline u}_0\over \p\tau}+{\overline u}_0\cdot{\nabla {\overline u}_0}
+ \nabla {\overline p}_1=
\pi (e_n{\overline \theta}_0),
\]
where $\tau=t/\sqrt\eps$ and $\nabla = \nabla_y$, with $y = x/\eps$.

\subsection{Problem setting}

We let $\om$ be a rectangular box, of thickness $h$
and Lebesgue measure $m(\om)$, in ${\bf R}^n$, $n=2$ or
$3$.
By $(e_i)$, $i=1,2$ or $i=1,2,3$, 
we denote the canonical basis in ${\bf R}^2$ or ${\bf R}^3$, respectively.
The system (1.2) is equipped with the following initial data:
\[
u(x,0)=u_0(x)\;\;{\rm and}\;\;T(x,0)=T_0(x),
\]
which are assumed to belong to $L^2(\om)$,
and boundary data:
\[
u=0\;\;{\rm at}\;\; x_n=0\;\;{\rm and}\,{\rm at}\;\;x_n=h.
\]
As above
\[T=T_1\;\;{\rm at}\;\;x_n=0\;\;{\rm and}\;\;T=T_2\;\;
{\rm at}\;\;x_n=h.
\]
Moreover we assume that
\[
u,\;\;\nabla u,\;\;T,\;\;\nabla T\;\;{\rm and}\;\;p
\]
are periodic
with period $l$ in the horizontal $x_1$-direction, in the two-dimensional case
and periodic with period $l$ in the horizontal $x_1$-direction
and  period $L$
in the $x_2$-direction in the three-dimensional case. In fact we will without
loss
of generality assume that $l=L$ throughout the paper.
This determines the pressure $p$ up to a constant that can be fixed
by normalization, see Remark 2.

Since we are interested in the fluctuations in the temperature we
follow \cite{fmt} and put
\[
\theta = (T-T_1-{x_n\over h}(T_2-T_1))
\]
and 
\[
p = p -g \alpha (x_n + {x_n^2 \over {2h}})(T_1-T_2).
\]

We get the
following system which is equivalent to (\ref{eq:bouss}). 
\begin{equation}
\label{eq:bouss1}
\left\{ \begin{array}{l}
\df {{\p{u}\over\p{t}}+(u\cdot{\nabla})u-\nu\Delta{u}
+\nabla{{p}}
= g\alpha e_n\theta,}\\[2ex]
\df {{\p{\theta}\over\p{t}}+(u\cdot{\nabla})\theta-
\kappa\Delta{\theta} = {{(T_1-T_2)} \over h} (u)_n},\\[2ex]
{\rm div}\,u=0,
\end{array} \right.\hfill x\in\om,\;t\in{\bf R}^+.\hfill 
\end{equation}
For the temperature
$\theta$ we get initial data $\theta(x,0)=\theta_0(x)$
and the new boundary data
$\theta=0$ at $x_n=0$ and
at $x_n=h$. The initial and boundary data for $u$
remain unchanged. Moreover, $u$ and
$\theta$ and their gradients and ${p}$ are periodic as above.
We will find it useful to work with the system
in the form (\ref{eq:bouss}) in some situations,
whereas the form (\ref{eq:bouss1}) is more suitable in other situations.

\section{The Navier-Stokes equation}

We will denote by $|\cdot|$ and $\|\cdot\|$ the usual norms in
$L^2(\om)$ and $H^1(\om)$. We denote by $\|\cdot\|_2$ the
$H^2(\om)$-norm and by $|\cdot|_{\infty}$ the $L^\infty(\om)$-norm. 
Further, $|u|_{2,\infty}={\rm ess}\,\sup|u(t)|$, where supremum is taken over all $t\geq 0$.

Let us consider the Navier-Stokes equation for incompressible fluids
\begin{equation}
\left\{ \begin{array}{l}
\df {{\p{u}\over\p{t}}+(u\cdot{\nabla})u
-\nu\Delta{u}+\nabla p = f,}\\[2ex]
{\rm div}\,u = 0, 
\end{array} \right.\;\;x\in\om,\;t\in{\bf R}^+,
\end{equation}
where $u$ is the velocity, $p$ is the pressure and $\nu$ is the viscosity,
with initial condition
\[
u(x,0)=u_0(x)
\]
and vanishing boundary conditions on $\p\om$ (periodic boundary conditions 
with mean zero, if $\om=T^3$). $f$ denotes the forcing and $\lambda_1$
is the smallest eigenvalue of $-\Delta$ on $\om$, with vanishing
boundary conditions on $\p\om$.
We start by an estimate which goes back to Leray \cite{lr1}. For the readers
convenience we present the (old) proof since the arguments therein
will be used repeatedly
both in the proof
of Theorem 2.2 and Theorem 3.3 below.
\begin{lem}
Every weak solution $u$ to the Navier-Stokes equation (2.1) satisfies
the estimate
\[
|u(t)|\leq|u_0| e^{-\lambda_1\nu t}+{|f|_{2,\infty}\over{\lambda_1\nu}}
(1-e^{-\lambda_1\nu t}).
\]
Moreover, there exists a sequence $t_j\to\infty$ such that
\[
\|u(t_j)\|^2\leq 3{|f|_{2,\infty}^2\over{\lambda^2_1\nu^3}}.
\]
\end{lem}
\begin{proof}
We take the inner product of (2.1) with $u$ and integrate over $\om$.
By the divergence theorem and the incompressibility, we are left with
\[
{1\over2}{d\over dt}|u(t)|^2+\nu|\nabla u(t)|^2=
\int_\om\, f(t)\cdot u(t)\, dx.
\]
The Schwarz and Poincar\'e inequalities give
\[
{1\over2}{d\over dt}|u(t)|^2+\lambda_1\nu|u(t)|^2\leq
|f(t)||u(t)|,
\]
so by cancellation of $|u(t)|$,
\[
{d\over dt}|u(t)|+\lambda_1\nu|u(t)|\leq
|f(t)|.
\]
An integration over $(0,t)$, taking sup in $t$ of $f$, gives
\begin{equation}
|u(t)|\leq|u_0| e^{-\lambda_1\nu t}+{|f|_{2,\infty}\over{\lambda_1\nu}}
(1-e^{-\lambda_1\nu t}).
\end{equation}
Next we integrate the inequality
\[
{1\over2}{d\over dt}|u(t)|^2+\nu|\nabla u(t)|^2
\leq
|f(t)||u(t)|
\]
over the interval $[t_1,t_2]$,
\[
\nu\int_{t_1}^{t_2}|\nabla u(t)|^2\, dt\leq
{1\over 2}(|u(t_1)|^2-|u(t_2)|^2)+
\int_{t_1}^{t_2}|f(t)||u(t)|\, dt.
\]
By (2.2) we get
\[
\nu\int_{t_1}^{t_2}|\nabla u(t)|^2\, dt\leq
|u_0|(\frac{|u_0|}{2}+|f|_{2,\infty})(e^{-\lambda_1\nu t_1}+e^{-\lambda_1\nu t_2})
+{|f|_{2,\infty}^2\over{(\lambda_1\nu)^2}}(1+(t_2-t_1)).
\]
Finally, we let $t_2=t_1+1$, and choose $t$ sufficiently large
to get
\[
\int_{t_1}^{t_1+1}\|u(t)\|^2\, dt\leq
{3|f|_{2,\infty}^2\over{\lambda_1^2\nu^3}}.
\]
This implies that there is a set of positive Lebesgue measure in every
interval $[t_1,t_1+1]$ such that
\[
\|u(t)\|^2\leq
{3|f|_{2,\infty}^2\over{\lambda_1^2\nu^3}},
\]
for $t$ in this set.
\end{proof}

We continue by stating a local existence theorem.
\begin{thm} Suppose that the pressure is normalized, i.e.,
\[
{1\over m(\om)}\int_\om\,p\,dx = 0,
\]
then there exists a unique local solution
$(u,\,p)$ of (2.1) in $ C([0,t];(H^1(\om))^{n+1})$, $n=2,\;3$, with initial data
$u_0=u(x,0)$ in $(H^1(\om))^{n}$. The local existence time $t$ depends
only on the $L^2$-norms of ${\rm curl}\,u_0$.
\end{thm} 
\noindent
Kreiss and Lorentz \cite{kl1} can be consulted for details of the proof.

It is well-known that global solutions of the Navier-Stokes equations exist,
\linebreak[4]
$ u \in C({\bf R^+};(H^1(\om))^{2})$, 
in two dimensions. In three dimensions
that analogous statement is open. However, if transients are allowed to settle
for a sufficiently long time, then global solutions exist in three
dimensions, after this settling of the initial velocity.

We now state and prove a result saying that every weak global solution to the
three-dimensional Navier-Stokes equation becomes a strong solution after some finite time
$t_0>0$. This also goes back to Leray \cite{lr1}.
The existence of a global attractor is due to Ladyshenskaya \cite{la1,la2}
and more recently Temam \cite{te1}, Sell \cite{Se96}
and P. L. Lions \cite{pll}. Both the statement and
the proof of Theorem 2.2 below are new and they are crucial in the statement
and proof of the new existence Theorem 3.3 below, for the Rayleigh-B\'enard
problem in the large aspect ratio.
\begin{thm} Consider the three-dimensional Navier-Stokes equation (2.1). 
For every weak solution $u_w\in C_w({\bf R}^+;(L^2(\om))^3)$ to (2.1) there 
exists a $t_0<\infty$, such that if the forcing $|f|_{2,\infty}$ is sufficiently
small there exists a unique strong solution 
$u\in C([t_0,\infty);(H^1(\om))^3)$ to (2.1), with data $u(t_0)=u_w(t_0)$. 
Moreover, (2.1) possesses a global attractor. If the initial data $|u_0|$ is 
small, then the solutions have global existence, i.e. $t_0=0$.
\end{thm}
\begin{proof}
The subscript $w$ indicates that the solutions are only weakly
continuos in $t$. We take the inner product of (2.1) with $\Delta u$ and integrate over $\om$.
Using Schwarz's inequality we obtain
\[
{1\over 2}{d\over dt}|\nabla u(t)|^2+\nu|\Delta u(t)|^2=
\int_\om\,f(t)\Delta u(t)\, dx-\int_\om((u(t)\cdot\nabla)u(t))\cdot\Delta u(t)\, dx
\]
\[
\qquad\leq|f(t)||\Delta u(t)|+|u(t)\cdot\nabla u(t)||\Delta u(t)|.
\]
Now
\[
|u(t)\cdot\nabla u(t)|\leq|u(t)|_\infty|\nabla u(t)|\leq
{K}|u(t)|^{1/4}\|u(t)\|_2^{3/4}|\nabla u(t)|
\]
by the Gagliardo-Nirenberg inequalities, where $K$ is a constant.
Thus,
\[
{1\over 2}{d\over dt}|\nabla u(t)|^2+\nu|\Delta u(t)|^2\leq
(|f(t)|+{K}|u(t)|^{1/4}\|u(t)\|_2^{3/4}|\nabla u(t)|)|\Delta u(t)|.
\]
An application of the inequality $ab\leq a^2/2\nu + \nu b^2/2$,
on the right hand side, gives
\[
{1\over 2}{d\over dt}|\nabla u(t)|^2+{\nu\over 2}|\Delta u(t)|^2\leq
{1\over\nu}(|f(t)|+{K}|u(t)|^{1/4}\|u(t)\|_2^{3/4}|\nabla u(t)|)^2
\]
\[
\qquad\leq \frac{1}{\nu}(|f(t)|^2+K^2|u(t)|^{1/2}\|u(t)\|_2^{3/2}|\nabla u(t)|^2)
\]
by another application of the inequality above with $\nu=1$.
From Lemma 2.1 we get, again by the same inequality,
\[
{1\over 2}{d\over dt}|u(t)|^2+\nu|\nabla u(t)|^2\leq
{1\over 2}(|f(t)|^2+|u(t)|^2).
\]
Adding these inequalities, applying the Sobolev inequality and repeating
the use of the inequality above results in
\[
{1\over 2}{d\over dt}\|u(t)\|^2+{\nu\over 4}\|u(t)\|_2^2\leq
{C_1}|u(t)|^2\|u(t)\|^8+C_2(|f(t)|^2+|u(t)|^2),
\]
where $C_1$ and $C_2$ are constants. Using Poincar\'e's inequality
and Lemma 2.1 we find that
\[
{1\over 2}{d\over dt}\|u(t)\|^2+{\lambda_1\nu\over 4}\|u(t)\|^2
-2{C_1}(|f|_{2,\infty}^2+|u_0|^2e^{-2\lambda_1\nu t})\|u(t)\|^8
\]
\[
\qquad\leq
3{C_2}(|f|_{2,\infty}^2+|u_0|^2 e^{-2\lambda_1\nu t}),
\]
since
\[
(|f|_{2,\infty}+|u_0|e^{-\lambda_1\nu t})^2\leq
 2(|f|_{2,\infty}^2+|u_0|^2e^{-2\lambda_1\nu t}).
\]
Now the point is that if the coefficient \, $|f|_{2,\infty}^2+|u_0|^2e^{-2\lambda_1\nu t}$ \,
is small, or the forcing is small and we have waited a sufficiently long time, to
let the initial data $|u_0|^2e^{-2\lambda_1\nu t}$ decay, then the inequality gives
us a bound on $\|u(t)\|$. The argument is as follows. We integrate the inequality
above over $[t_0,t]$ to get
\[
\|u(t)\|^2-4C_1
\int_{t_0}^t (|f|^2_{2,\infty}+|u_0|^2e^{-2\lambda_1\nu t_0})
e^{-\beta(t-s)}\|u(s)\|^8\, ds
\]
\[
\qquad\leq
{6{C_2}\over\beta}(|f|_{2,\infty}^2+|u_0|^2 e^{-2\lambda_1\nu t_0})(1-e^{-\beta(t-t_0)}),
\]
where $\beta=\lambda_1\nu/2$.
Now assume that $\|u(s)\|$ assumes its maximum in the interval $[t_0,t]$ at $s=t$.
Then
\[
\|u(t)\|^2-\frac{4C_1}{\beta}
(|f|^2_{2,\infty}+|u_0|^2e^{-2\lambda_1\nu t_0})
\|u(t)\|^8\, ds
\]
\[
\qquad\leq
{6{C_2}\over\beta}(|f|_{2,\infty}^2+|u_0|^2 e^{-2\lambda_1\nu t_0})(1-e^{-\beta(t-t_0)}).
\]
Now put
\[
v(t)=\|u(t)\|^2, \;\;a=\frac{4C_1}{\beta}
(|f|^2_{2,\infty}+|u_0|^2e^{-2\lambda_1\nu t_0})
\]
and
\[
M={6{C_2}\over\beta}(|f|_{2,\infty}^2+|u_0|^2 e^{-2\lambda_1\nu t_0}).
\]
Then the inequality above can be written as
\[
v(t)-av^4(t)\leq M.
\]
\begin{figure}[htbp]
\centering 
\includegraphics[height=3in]{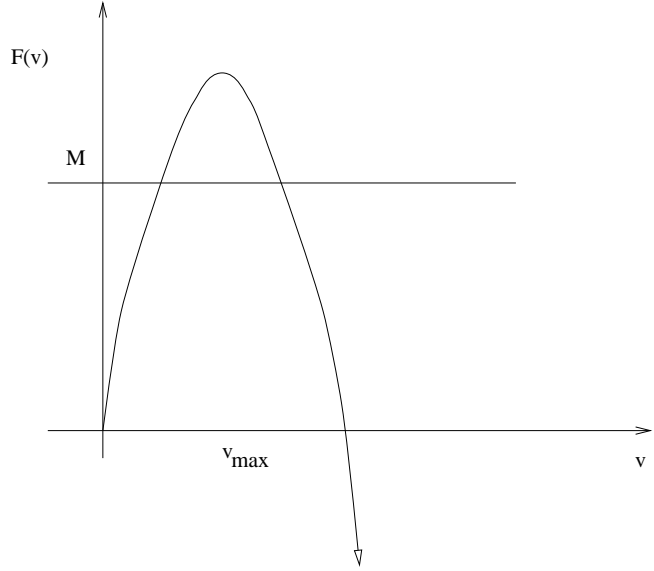}
\caption{The bound on the norm.}
\label{fig:bound}
\end{figure}
The graph of the function $F(v)=v-av^4$ is shown in Figure 1. It is concave and attains
its maximum at $v_{max}=1/(4a)^{1/3}$. By Lemma 2.1 there exists a $t_0 > 0$ such that
\[
v(t_0)\leq\|u(t_0)\|^2\leq\frac{3|f|_{2,\infty}}{\lambda_1^2\nu^3}
\]
and we can choose $|f|_{2,\infty}$ so small that $v(t_0)$ lies between
zero and $v_{max}$, see Figure 1. Moreover, $v(t)$ can never reach $v_{max}$,
because
\[
F(v(t))=v(t)-av^4(t)\leq M < \frac{3}{4}\frac{1}{(4a)^{1/3}}=F(v_{max}).
\]
It is clear that we can choose $t_0$ so large that
\[
M={6{C_2}\over\beta}(|f|_{2,\infty}^2+|u_0|^2 e^{-2\lambda_1\nu t_0})
\]
\[
\qquad\leq
\frac{3}{4}\frac{1}{(4a)^{1/3}}=\frac{3}{4}
\frac{1}{(\frac{16C_1}{\beta}
|f|^2_{2,\infty}+|u_0|^2e^{-2\lambda_1\nu t_0})^{1/3}}.
\]
Moreover, notice that the derivative of $F$ in Figure 1 is positive at the point, where
$F$ first reaches $M$, i.e.,
\[
F'(v)=1-4av^3=1+\frac{4(v-av^4)}{v}-4=\frac{4M}{v}-3\geq 0,
\]
which gives the bound 
\[
v(t)\leq \frac{4M}{3}+\delta
\]
where $\delta$ is arbitrarily small. The only question left is whether
the initial data makes sense as a function in $H^1(\om)$. But we have
already used above that, by Lemma 2.1, there exists a sequence
$t_j\to\infty$ such that $\|u(t_j)\| < \infty$. We now choose the
initial time to be the smallest $t_j\geq t_0$. Then we let this
$t_j$ be the new $t_0$.

 Now we apply the local existence Theorem 2.1 and the above bound to 
get global existence. We recall the definition of an absorbing set from Coddington and Levinson 
\cite{CL55}. A
set $D \subset L^2$ is an absorbing set if for every bounded set $M$ there exists a time $t(M)$ such that 
$t > t(M)$ implies that $u(t) \in D$, if $u(0) \in M$. The estimate in Lemma 2.1 shows that the 
weak-flow of the weak solution of the Navier-Stokes equations has an absorbing set in $L^2$. But we have shown
in addition that there exists a $t_1(M)$ such that the flow is a continuos flow of strong solutions
for $t > t_1(M)$ and has
the absorbing set $D \subset H^1$. 
Moreover, since $H^1(\om)$ is compactly
embedded in $L^2(\om)$ it follows that the Navier-Stokes equation has a global
attractor in $L^2(\om)$. In fact one can now show, see \cite{kl1}, 
that the solutions are spatially
smooth and thus $D \subset C^\infty$. It follows, see Hale \cite{key8}, Babin and Vishik \cite{key1}, Temam 
\cite{te1} and Birnir and Grauer \cite{bg2}, that the Navier-Stokes equation has a global attractor 
consisting of spatially smooth solutions. 
It is also clear that if $|u_0|$ is small then $v(t)$
satisfies the above bound and the solutions have global existence for $t_0=0$.
\end{proof}

\begin{rem}
The statement of Theorem 2.2 is really just a restatement of 
Leray's classical result of global existence for small initial data.
The duality
in the smallness condition is that one can either take small initial data or,
with arbitrary initial data, just wait for a sufficiently long time. 
This statement is what
most physicist and engineers are interested
in, at least for low Reynolds numbers. Namely, transients die out and the flow settles down to the flow 
on the attractor in at most a few hours, in experiments, see for instance
\cite{hu2}.
\end{rem} 
\begin{rem}
The existence of the pressure term follows from a standard orthogonality
argument which we present below. 
Consider the Navier-Stokes equation (\ref{eq:bouss})
and the solution given by Theorem 2.2. We take the inner product by
$u$, integrate over $\om$, apply the divergence theorem, use the incompressibility
and collect all terms on one side. This gives
\[
\int_\om(f-{\p u\over \p t}-u\cdot\nabla u +\nu\Delta u)\cdot u\,dx=0.
\]
Now we recall that the orthogonal to the divergence free elements in $L^2$
are gradients in $L^2$. Therefore, there exists a unique gradient $\nabla p$
in $L^2(\om)$ given by
\[
\nabla p=f-{\p u\over \p t}-u\cdot\nabla u+\nu\Delta u.
\]
\end{rem}

We say that a set $M$ is {\em invariant} under the flow, defined by a nonlinear
semi-group $S(t)$, if $S(t)M=M$,
i.e. $M$ is both positively and negatively invariant.
An {\em attractor} ${\cal A}$ is an invariant
set which attracts a neighbourhood $U$ or for all $x_o\in U$,
$S(t)x_o$ converges to ${\cal A}$ as $t \rightarrow \infty$.  The largest
such set $U$ is called the {\em basin of attraction} of the attractor
${\cal A}$.  If the basin of attraction of an attractor ${\cal A}$ 
contains all bounded sets of $X$ and ${\cal A}$ is compact, 
then we say that ${\cal A}$ is the {\em global attractor}.

We have shown above that the Navier-Stokes equations have a global attractor.
But we have not explained what the semi-group is that shrinks bounded sets onto
the attractor. It is possible to define a semi-group on the space of weak-solutions,
see Sell \cite{Se96}, but in light of Theorem 2.2 we can do much better. 
Namely, a bounded subset of a Banach space is a complete metric space and we have shown that for
every bounded set $M \subset L^2$ eventually the solution starting in $M$ will lie in an
absorbing set $D \subset C^\infty$. It is the semi-group defined on this complete metric space
$D \subset H^1$ that has the global attractor
${\cal A} = \omega (D)$, where $\omega$ denotes the $\omega$-limit set, see \cite{bg2}.
$\cal A$ is invariant, so on it we can solve the Navier-Stokes equation in backward time, it
is non-empty and compact and 
attracts a neighbourhood of itself, see \cite{bg2}, moreover, it also attracts every bounded set in $L^2$.
This attractor
consists of spatially smooth solutions and has finite Hausdorff
and fractal dimensions, see \cite{l1, ma, r1, r2}.  The dimension estimates are however large,
see \cite{te1}, 
and not necessarily a good indication of the size of the attractor.

We focus on the core of the attractor ${\cal A}$ which is called
the basic attractor ${\cal B}$.  
For now, let us assume that the Banach space $X$ is finite dimensional.
An attractor ${\cal B}$ is called a {\em basic attractor} if
\begin{itemize}
\item{The basin of attraction of ${\cal B}$ has positive measure.}
\item{There exists no strictly smaller ${\cal B}' \subset {\cal B}$,
such that up to sets of measure zero, {\em basin}(${\cal B}$)
$\subset${\em basin}(${\cal B}'$).}
\end{itemize}
A {\em global basic attractor} is a basic attractor whose basin of 
attraction contains all bounded sets of $X$ up to sets of measure
zero.

A theorem by Milnor \cite{ml} states that in finite dimensions
there exists a unique decomposition of the attractor ${\cal A}$,
\[
        {\cal A}={\cal B}\cup{\cal C}
\]
where ${\cal B}$ is a basic attractor and ${\cal C}$ is a remainder
such that $m(basin({\cal C}) \backslash basin({\cal B}))=0$
where $m$ is the Lebesque measure.  

The notion of a basic attractor can be extended to infinite dimensions
and Milnor's Theorem was first proven in an infinite-dimensional setting
by Birnir and Grauer \cite{bg1}. They used 
cumbersome projections to a finite dimensional space where
the $ {\cal A}$-attractor resides. A more elegant notion of 
a ${\cal B}$-attractor
requires an extension of the concepts of measure zero and almost everywhere.
Their counterparts in infinite dimensions are {\em shy} and {\it prevalent}
sets respectively, see Hunt {\it et al.} \cite{shy}.  These are
defined in the following way.
Let $X$ denote a Banach space.  We denote by $S+v$ the 
translate of the set $S\subset X$ by a vector $v\in X$.  
A measure $\mu$ is said to be {\em transverse} to a Borel set $S \subset X$
if the following two conditions hold:
\begin{itemize}
\item{There exists a compact set $U\subset X$ for which $0<\mu(U)<\infty$.}
\item{$\mu(S+v)=0$ for every $v\in X$.}
\end{itemize}
A Borel set $S\subset X$ is called {\em shy} if there exists a compactly
supported
measure transverse to $S$.  More generally, a subset of $X$ is called shy
if it is contained in a shy Borel set.  The complement of a shy set is
called a {\em prevalent} set.

The infinite-dimensional analogue of Milnor's Theorem can now be stated,
\begin{thm}
Let ${\cal A}$ be the compact attractor of a continuous map $S(t)$ on
a separable Banach space $X$.  Then ${\cal A}$ can be decomposed into
a maximal basic attractor ${\cal B}$ and a remainder ${\cal C}$,
\[
        {\cal A}={\cal B}\cup{\cal C}
\]
such that the realm of attraction of ${\cal B}$ is prevalent but the
realm of attraction of ${\cal C}$, excluding points that are attracted
to ${\cal B}$, is shy.
\end{thm}
\noindent
For a proof of this Theorem see Birnir \cite{bb2} and \cite{b4}. It implies
that the Navier-Stokes equation has a maximal  ${\cal B}$-attractor in 
two and three dimensions. If  ${\cal B}$ can be decomposed into finitely many 
(disjoint) minimal ${\cal B}$-attractors, then the union of the realms
of attraction of these minimal ${\cal B}$-attractors is the whole space 
$(L^2({ \om \subset \bf R}))^n, n = 2, 3.$
The realm of
attraction is a slight generalization of the basin of attraction, as the 
basin is an open set, but the realm can be either open or closed.

For systems that are simple enough, the basic attractor contains only
the stable trajectories of the global attractor.
When transients are ignored, one will only see the basic
attractor in physical  experiments and numerical simulations, not the
remainder ${\cal C}={\cal A}\backslash {\cal B}$ of the 
attractor.  All the relevant dynamics are therefore contained in the
basic attractor. Ladyshenskaya \cite{la3} gives more examples of ${\cal B}$-attractors
for the Navier-Stokes equation with nonlinear viscosity. 

\section{The Rayleigh-B\'enard convection}

\subsection{Existence results}

Now consider the Boussinesq system (1.1)
equipped with initial data
\[
u(x,0)=u_0(x)\;\;{\rm and}\;\;T(x,0)=T_0(x).
\]
The boundary conditions are periodic on the vertical surfaces
$x_k=0,\,l$ for $k<n$. Furthermore, on the horizontal surfaces
$x_n=0$ and $x_n=h$, the velocity $u$ vanishes (no-slip condition)
and the temperature $T$ is kept constant, i.e.,
\[
T=T_1\;\;{\rm on}\;\;x_n=0,\;\;T=T_2\;\;{\rm on}\;\;x_n=h.
\]

We want to state a global existence theorem for
the solution to the system (1.1). For that
purpose we start with the following maximum principle,
c.f. \cite{fmt}:
\begin{lem}
Suppose that $u$ and $T$ solve (1.1). If 
\[
T_2\leq T(x,0)\leq T_1,
\]
for a.e. $x\in\om$, then
\[
T_2\leq T(x,t)\leq T_1,
\]
for a.e. $x\in\om$ and all $t\geq 0$.
\end{lem}
\begin{proof}
We consider
\[
(T-T_1)_+(x,t) = {\rm ess}\,\sup_{x\in \om}(T-T_1)(x,t).
\]
The second equation in (3.1) gives
\[
\df {{\p\over\p{t}}{(T-T_1)_+}+(u\cdot{\nabla})(T-T_1)_+
-\kappa\Delta{(T-T_1)_+} = 0.}
\]
We multiply
this equation by $(T-T_1)_+$, integrate over $\om$ and use the divergence theorem,
to get
\[
{1\over 2}{d\over dt}|(T-T_1)_+(t)|^2+
\kappa|\nabla(T-T_1)_+(t)|^2=0.
\]
Thus, by  Poincar\'e's inequality
\[
{1\over 2}{d\over dt}|(T-T_1)_+(t)|^2+
\lambda_1\kappa|(T-T_1)_+(t)|^2\leq 0,
\]
where again $\lambda_1$ is the first eigenvalue of the negative
Laplacian with vanishing boundary conditions on $\Omega$. 
Consequently
\[
|(T-T_1)_+(t)|\leq
|(T-T_1)_+(0)| e^{-\lambda_1\kappa t},
\]
which shows that
\[
(T-T_1)_+(\cdot,t)=0,
\]
for all $t\geq 0$, if
\[
(T-T_1)_+(\cdot,0)=0.
\]
Similarly we get that
\[
(T-T_2)_-(x,t)={\rm ess}\,\sup_{x\in\om}(-(T-T_2)(x,t)) = 0,
\]
if $(T-T_2)_-(x,0)=0$ and we conclude that
\[
T_2\leq T(x,t)\leq T_1. 
\]
\vspace{-1.5em}
\end{proof}
\begin{rem}
Lemma 3.1 actually yields a uniform bound on $T$
in $L^2([0,s]\times\om)$ for any $s> 0$. This also gives a bound for 
$\theta$, see below.
\end{rem}
\begin{lem}
Every weak solution $u$ to the Boussinesq equations (1.4) satisfies
the estimate
\[
|u(t)|\leq|u_0| e^{-\lambda_1\nu t}+{{K h^{1/2}}\over{\lambda_1\nu}}
(1-e^{-\lambda_1\nu t}),
\]
\[
|\theta(t)|\leq {{K h^{1/2}}\over{g \alpha}},
\]
where $K = g\alpha (T_1-T_2)L/3^{1/2}$.
The equations possess an absorbing set in $(L^2(\Omega))^4$, defined by
\[
|u(t)|+|\theta(t)|\leq 
\left( \frac{1}{g\alpha} +\frac{1}{\lambda_1\nu} \right) K h^{1/2} + \delta,
\]
where $\delta$ is arbitrarily small.
Moreover, there exists a sequence $t_j\to\infty$ such that
\[
\|u(t_j)\|^2+\|\theta(t_j)\|^2\leq K_1,
\]
where 
\[
K_1 = 3{K^2\over{\lambda_1^2 \nu^2}}
\left( \frac{h}{\nu}+{{(T_1-T_2)^2}\over{\lambda_1^2 \kappa^3 h}} \right).
\]
\end{lem}
\begin{proof}
We recall the relationship between $T$ and $\theta$
\[
T = \theta + T_1 - {x_n \over h}(T_1-T_2).
\]
The maximum principle in Lemma 3.1, for T, $ T_2 \le T \le T_1$, implies
that 
\[
-(1-{x_n \over h})(T_1-T_2) \le \theta \le {x_n \over h}(T_1-T_2).
\]
Thus 
\[
|\theta|^2_2 \le (T_1-T_2)^2 L^2 h/3.
\]
Now consider the system (1.4).
We multiply the first equation
by $u$ and integrate over $\Omega$. Integration by parts and Schwarz's 
inequality give 
\[
{1\over2}{d\over dt}|u(t)|^2+\lambda_1\nu|u(t)|^2\leq
g \alpha |\theta (t)||u_n(t)|.
\]
Thus by the use of the bound for $\theta$ above, and Poincar\'e's inequality,
we get
\[
{d\over dt}|u(t)|+\lambda_1\nu|u(t)|\leq
g \alpha |\theta (t)| \le g \alpha (T_1-T_2) L (h/3)^{1/2}.
\]
Integration in t then gives
\[
|u(t)|\leq|u_0| e^{-\lambda_1\nu t}+ {{K h^{1/2}}\over{\lambda_1 \nu}}
(1-e^{-\lambda_1\nu t}),
\]
where  $K = g \alpha (T_1-T_2)^{1/2} L / (3)^{1/2}$. 

We combine the bounds for $\theta$ and $u$ to get the absorbing set
\[
|\theta| + |u| \le K h^{1/2}\left( \frac{1}{g\alpha} +\frac{1}{\lambda_1\nu}
\right)+\delta
\]
in $(L^2(\Omega))^4$, where
\[
K = {g\alpha}{{(T_1-T_2)} \over {3^{1/2}}} L  
\]
and $\delta$ is arbitrarily small.

The last statement of the lemma is proven by a straightforward application
of Lemma 2.1, to the equations (1.4), if we recall that the nonlinear term did not
play a role in the proof of Lemma 2.1. Namely, for the first equation in (1.4),
\[
|f|_{2,\infty}^2 = g^2 \alpha^2 |\theta|^2 = K^2 h,
\]
and for the second equation
\[
|f|_{2,\infty}^2 = {{(T_1-T_2)^2}\over {h^2}}|u|^2=
 {{(T_1-T_2)^2 K^2}\over {\lambda_1^2 \nu h}}.
\]
We divide these by the decay coefficients $\lambda_1^2 \nu^3$ and 
$\lambda_1^2 \kappa^3$, respectively and add them. This produces the bound
on the $H^1$ norm for the sequence ${t_j}$. 
\end{proof}

\begin{thm} Suppose that the pressure is normalized, i.e.,
\[
{1\over m(\om)}\int_\om\,p\,dx = 0,
\]
then there exists a unique local solution
$(u,\theta ,p)$ of (1.4) in $ C([0,t];(H^1(\om))^{n+2})$, $n=2,\;3$, with initial data
$(u_0, \theta_0)(x)$ in $(H^1(\om))^{n+1}$.
\end{thm}
\begin{proof}
The maximum principle in Lemma 3.1, for T, $ T_2 \le T \le T_1$, 
implies a maximum principle for $\theta$, 
\[
-(1-{x_n \over h})(T_1-T_2) \le \theta \le {x_n \over h}(T_1-T_2),
\]
by the relationship between $T$ and $\theta$
\[
T = \theta + T_1 - {x_n \over h}(T_1-T_2).
\]
This implies that
\[
|\theta|^2_2 \le (T_1-T_2)^2 L^2 h/3.
\]
The rest of the proof, using this bound on $\theta$, is similar to 
the proof of the local existence of the Navier Stokes equation. 
Kreiss and Lorentz \cite{kl1} can be consulted for details.
Then given the local solution $(u,\theta)(x,t)$ the pressure is recovered as 
in Remark 2. 
\end{proof}
In two dimensions we have the following global existence result.
\begin{thm}
The Boussinesq system (1.4) has a unique global solution
$(u,\,\theta)$\ in $\\ C({\bf R}^+;(H^1(\om))^3)$, where
$\om$ is a bounded open set in ${\bf R}^2$. Moreover, the system (1.4)
possesses a global attractor in $(L^2(\om))^3$.
\end{thm}
\begin{proof}
A proof of Theorem 3.2 can be found in Foias et. al. \cite{fmt}.
\renewcommand{\qedsymbol}{}
\end{proof}

In three dimensions the following global existence result holds true:
\begin{thm}
Consider the Boussinesq system (1.4). For every weak solution
$u_{w},\,\theta_{w}$ in
$C_w({\bf R}^+;(L^2(\om))^4)$, $\om$ a bounded open set in ${\bf R}^3$ with a large
aspect ratio, $ {L \over h} >> {g\alpha^2 (T_1-T_2)^2 L^3\over \nu\kappa} $,
where $L$ is the width (radius)
and $h$ is the height of $\Omega$;
there exist a time $t_0$, such that there exists a unique strong
solution $(u,\,\theta)$ in
$C([t_0,\infty);(H^1(\om))^4)$, $t\geq t_0$, with initial data
$u(t_0)=u_{w}(t_0)$ and $\theta(t_0)=\theta_{w}(t_0)$.
Moreover, the system (1.4) possesses a global attractor in $(L^2(\om))^4$.
\end{thm}
\begin{proof}
The subscript $w$ indicates that the weak solutions are only
weakly continuos in $t$. We multiply the second equation in (1.4) above by $\theta$ and treat it in the same 
way as $u$ in the proof of Lemma 2.1 to get

\[
{1\over2}{d\over dt}|\theta (t)|^2+\kappa|\nabla \theta (t)|^2\leq
{{(T_1-T_2)} \over h}|\theta (t)||u_n(t)|,
\]
and by adding and subtracting $\beta |\theta|$ and applying Poincar\'e's 
inequality, we get

\[
{d\over dt}|\theta (t)|+\beta |\theta (t)|+(\lambda_1^{1/2}\kappa
-{\beta \over {\lambda_1^{1/2}}}){|\nabla \theta (t)|\leq
{{(T_1-T_2)} \over h}}|u_n(t)|
\]
\[
\qquad\le {{(T_1-T_2)} \over h}(|u_0| e^{-\lambda_1\nu t}+ {{K
h^{1/2}}\over{\lambda_1\nu}}
(1-e^{-\lambda_1\nu t})) \le c_1 + c_2e^{-\lambda_1\nu t}.
\]

Then integrating with respect to t, 
we get
\[
(\lambda_1^{1/2}\kappa
-{\beta \over {\lambda_1^{1/2}}})\int_0^t e^{-\beta (t-s)}|\nabla \theta (s)| ds
\]
\[
\qquad\le |\theta (0)|e^{-\beta t} + {c_1 \over \beta} (1-e^{-\beta t}) 
+c_2 {{(e^{-\beta t} - e^{-\lambda_1 \kappa t})} 
\over {(\lambda_1 \kappa - \beta)}}.
\]

Next we multiply the u equation in (1.4) by $\Delta u$ and integrate over $\Omega$.
By integration by parts and Schwarz's inequality

\[
{1\over 2}{d\over dt}|\nabla u(t)|^2+\nu|\Delta u(t)|^2 
\leq g\alpha |\nabla \theta (t)||\nabla u_n (t)|+|u(t)\cdot\nabla u(t)|
|\Delta u(t)|
\]
\[
\qquad\le g\alpha |\nabla \theta (t)||\nabla u_n (t)| + 
C(1+\lambda_1^{-1}+\lambda_1^{-2})^{3/8}| u(t)|^{1/4}|\nabla u (t)| 
|\Delta u (t)|^{7/4},
\]
by the Gagliardo-Nirenberg inequalities, where we have used that the 
$H^2$ Sobolev norm is bounded by
\[
\|u\|_2 \le C(| u(t)|^2+|\nabla u (t)|^2 +
|\Delta u (t)|^2)^{1/2} \le C(1+\lambda_1^{-1}+\lambda_1^{-2})^{1/2}|\Delta u (t)|,
\]
by Poincar\'e's inequality. We use Young's inequality to eliminate 
$|\Delta u|^{7/4}$, namely
\[
(({\nu \over 2})^7 C_0)^{1/8}|u|^{1/4}|\nabla u| |\Delta u |^{7/4} \le 
C_0 |u|^2|\nabla u|^8+ {\nu \over 2}|\Delta u |^2
\]
so 
\[
{1\over 2}{d\over dt}|\nabla u(t)|^2+{\nu \over 2}|\Delta u(t)|^2 
\leq g\alpha |\nabla \theta (t)||\nabla u (t)|+ C_0 |u(t)|^2|\nabla u(t)|^8
\]
and by Poincar\'e's inequality
\[
{d\over dt}|\nabla u(t)|+{{\lambda_1  \nu} \over 2} |\nabla u(t)| 
\leq  g\alpha |\nabla \theta (t)|+
C_0 |u(t)|^2|\nabla u(t)|^7.
\]
Then we integrate the equation
\[
{d\over dt}|\nabla u(t)|+\beta|\nabla u(t)|
-C_0 |u(t)|^2|\nabla u(t)|^7 \leq  g\alpha |\nabla \theta (t)|,
\]
where $\beta = \lambda_1 \nu/2$. 
We integrate from the initial time $t_0$, to get

\[
|\nabla u(t)|-C_0 \int_{t_0}^t 
|u(s)|^2|\nabla u(s)|^7e^{-\beta (t-s)}ds 
\]
\[ 
\qquad\le c_1 + c_2e^{-\beta (t-t_0)} + c_3 e^{-\lambda_1 \nu (t-t_0)},
\]
by the above inequality for $\int |\nabla \theta| dt$.
Now if $|\nabla u(s)|$ attains its maximum on the interval $[t_0,t]$ at $s=t$,
then
\[
\int_{t_0}^t 
|u(s)|^2|\nabla u(s)|^7e^{-\beta (t-s)}ds \leq
|\nabla u(t)|^7\int_{t_0}^t 
|u(s)|^2e^{-\beta (t-s)}ds 
\]
\[
\qquad\leq\frac{|\nabla u(t)|^7}{2}\int_{t_0}^t 
\left(|u_0|^2e^{-3\beta s}e^{-\beta t}+ {{K^2h}\over{\lambda_1^2 \nu^2}}(1-e^{-2\beta s})^2e^{-\beta (t-s)}\right)ds
\]
\[
\qquad\leq \frac{|\nabla u(t)|^7}{6\beta}\left(|u_0|^2(e^{-3\beta t_0}-e^{-3\beta t})e^{-\beta t}+
3{{K^2h}\over{\lambda_1^2 \nu^2}}(1-e^{-\beta (t-t_0)})\right),
\]
where we have used the bound
\[
|u(t)|\leq|u_0| e^{-\lambda_1\nu t}+ {{Kh^{1/2}}\over{\lambda_1 \nu}}
(1-e^{-\lambda_1\nu t}),
\]
from Lemma 3.2 and the inequality $ab\leq a^2/2+b^2/2$.
Thus, if we put $|\nabla u(t)|=v(t)$, we obtain the inequality
\[
v(t)-av^7(t)\leq M.
\]
The constants $a$ and $M$ are,
\[
a= {{K^2h}\over{2 \lambda_1^2 \nu^2 \beta}}= {g^2\alpha^2}{{(T_1-T_2)^2} \over
{6 \lambda_1^2 \nu^2 \beta}} L^2 h
\]
and 
\[
M = 2c_1 = 2{{(T_1-T_2)K}\over{h^{1/2}\lambda_1\nu}}=
2g\alpha (T_1-T_2)^{2}L/3^{1/2}\lambda_1\nu h^{1/2}.
\]
This means that for a large aspect ratio $L/h > >
{g\alpha^2(T_1-T_2)^2 L^3\over \nu\kappa}$, $a$ 
becomes small.
Now we repeat the arguments from the proof of Theorem 2.2 and conclude that
the function $F(v)=v-av^7$ is concave and attains its maximum at 
$v_{max}=1/(7a)^{1/6}$.
This maximum is
\[
F(v_{max})=\frac{6}{7}\frac{1}{(7a)^{1/6}},
\]
so that $v(t)$ can not escape beyond its maximum, see Figure 1.
Moreover arguing as in Theorem 2.2 we conclude that the derivative of $F$ is positive at
the point where $F$ first reaches $M$ and this gives us the bound
\[
v\leq \frac{7M}{6}.
\]

The last step is to get a bound for $|\nabla \theta|$. We multiply the 
$\theta $ equation in (1.4) by $\Delta \theta$, to get
\[
{1\over 2}{d\over dt}|\nabla \theta(t)|^2+\kappa | \Delta \theta (t)|^2 
\leq {{(T_1-T_2)} \over h} |\nabla \theta (t)||\nabla u (t)|
+|u(t)\cdot\nabla \theta (t)|
|\Delta \theta(t)|
\]
by Schwarz's inequality
\[
\qquad\le {{(T_1-T_2)} \over h} |\nabla \theta (t)||\nabla u (t)| 
+ |u|_6|\nabla \theta (t)|_3|\Delta \theta(t)|, 
\]
by H\"{o}lder's inequality
\[
\qquad\le  {{(T_1-T_2)} \over h} |\nabla \theta (t)||\nabla u (t)| + 
C(|\nabla u (t)|^2|\theta|/\delta^3 +\delta |\Delta \theta|)
 |\Delta \theta (t)|,
\]
by Poincar\'e's and Sobolev's inequalities, because 
\[
|u|_6 \le K \|u\| \le C |\nabla u|, 
\] 
as above, and because,
\[
|\nabla \theta|_3 \le C \|\nabla \theta\|_{1/2} \le K ({|\nabla \theta|
\over \delta} + \delta |\Delta \theta|) \le C ({|\theta|
\over {\delta^3}} + \delta |\Delta \theta|)
\]
by two applications of interpolation, where $\delta$ is small. Thus
\[
{1\over 2}{d\over dt}|\nabla \theta(t)|^2+\kappa | \Delta \theta (t)|^2 
\leq {{(T_1-T_2)} \over h} |\nabla \theta (t)||\nabla u (t)| + 
{{C^2|\nabla u (t)|^4|\theta|^2} \over {\delta^6}} +2\delta |\Delta \theta (t)|^2,
\]
by Young's inequality. Now moving the $|\Delta \theta (t)|^2$ term over 
to the left hand side of the inequality and applying Poincar\'e's inequality,
we get
\[
{1\over 2}{d\over dt}|\nabla \theta(t)|^2+\lambda_1(\kappa-3\delta)|\nabla \theta(t)|^2 
\leq {1\over 2}{d\over dt}|\nabla \theta(t)|^2+
(\kappa-3\delta)|\Delta \theta(t)|^2 \le c_3,
\]
since
\[
{{(T_1-T_2)} \over h} |\nabla \theta (t)||\nabla u (t)| \leq \lambda_1 \delta
|\nabla \theta (t)|^2 + {{(T_1-T_2)^2} \over {4 h^2 \lambda_1 \delta}}
|\nabla u (t)|^2,
\]
by the inequality $ab\leq a^2+b^2/4$.
Thus 
\[
|\nabla \theta(t)|^2 \le |\nabla \theta(t_0)|^2e^{-\gamma (t-t_0)} 
+ (2c_3/\gamma)(e^{-\gamma t_0}-e^{-\gamma t}),
\]
where $\gamma = 2(\kappa -3\delta)$.

We can now put the estimates for $|\nabla u|$ and $|\nabla \theta|$ together
to get an absorbing set
\[
|\nabla u (t)| + |\nabla \theta (t)| \le constant +\eps,
\]
where $\eps$ is arbitrarily small for $t$ large enough. 
Namely, by Lemma 3.2, there exists a $t_j$, such that 
$(u,\theta)(t_j) \in H^1(\om)$ and the global bound on the $H^1$ norm holds
for $t=t_0=t_j$. Combined with the local existence Theorem 3.1, the a priori 
bound 
above now gives the existence of a global solution and an absorbing set 
in $H^1(\om)$, for $t \ge t_0$.
The existence of
a global attractor in $(L^2(\Omega))^4$ then follows from the compact embedding
of $(H^1(\Omega))^4$ in $(L^2(\Omega))^4$,  see Hale \cite{key8}, Babin and Vishik \cite{key1}, Temam 
\cite{te1} and Birnir and Grauer \cite{bg2}. 
\end{proof}
\begin{cor} 
The Boussinesq system (1.4) has a maximal ${\cal B}$-attractor whose
basin is a prevalent set in
$(L^2({\bf R}^n)^{n+1})$, $n=2,\,3$.
\end{cor}
The proof is a straight-forward application of Theorem 2.3

Theorem 3.3 says that in experiments in pattern formation where one has a 
large aspect ratio, it is only necessary to wait a short time to have global
spatially smooth solutions. The smoothness follows from the smoothness of the 
nonlinear semigroup, see Kreiss and Lorentz \cite{kl1}. Corollary 3.1 says that 
what is observed in the experiments after the initial settling-down time, is 
a $\cal B$-attractor that is a component of the maximal $\cal B$-attractor. In 
other word this attractor is observed for an open set of initial conditions. 
(However, not necessarily an open set of parameter values.) It is then useful
to know how long one has to wait and this time is easily estimated. 
It is the time it takes the exponentially decaying initial data to become
comparable in size with the $L^2$ absorbing set. Namely,
\[
t= \frac{1}{\lambda_1 \nu} \ln \left[ {{\lambda_1 \nu |u_0|}\over{Kh^{1/2}}} \right]
= \frac{1}{\lambda_1 \nu} \ln \left[ {{3^{1/2} \lambda_1 \nu |u_0|}\over{
g \alpha (T_1 - T_2) L h^{1/2}}} \right].
\]
In experiments for Prandtl numbers close to one, this time is measured in hours
for experiments that take days \cite{hu2}.

\section{Scaling and expansions}

As a starting point we introduce a small "scaling" parameter $\eps>0$ and consider
the scaled system 
\begin{equation}
\left\{ \begin{array}{l}
\df {{\p{u_\eps}\over\p{t}}+(u_\eps\cdot{\nabla})u_\eps
-\eps^{\gamma} \nu\Delta{u_\eps}+\nabla p_\eps = g\alpha(T_\eps-T_2),}\\[2ex]
\df {{\p{T_\eps}\over\p{t}}+(u_\eps\cdot{\nabla})T_\eps
-\eps^{\gamma} \kappa\Delta{T_\eps} = 0,}\\[2ex]
{\rm div}\,u_\eps = 0, 
\end{array} \right.\;\;x\in\om,\;t\in{\bf R}^+,
\end{equation}
equipped with initial data
\[
u_\eps(x,0)=a(x)\;\;{\rm and}\;\;T_\eps(x,0)=b(x).
\]
The boundary conditions are periodic on the vertical surfaces
$x_k=0,\,l$ for $k<n$. Furthermore, on the horizontal surfaces
$x_n=0$ and $x_n=h$, the velocity $u_\eps$ vanishes (no-slip condition)
and the temperature $T_\eps$ is kept constant, i.e.,
\[
T_\eps=T_1\;\;{\rm on}\;\;x_n=0,\;\;T_\eps=T_2\;\;{\rm on}\;\;x_n=h.
\]
A small value of $\eps$ corresponds to a high 
Reynolds number
and the solutions will likely become turbulent.
We also consider the scaled system corresponding to (1.4), i.e.,
\begin{equation}
\left\{ \begin{array}{l}
\df {{\p{u_\eps}\over\p{t}}+(u_\eps\cdot{\nabla})u_\eps-\eps^{\gamma}\nu\Delta{u_\eps}
+\nabla{{p_\eps}}
= g\alpha e_n\theta_\eps,}\\[2ex]
\df {{\p{\theta_\eps}\over\p{t}}+(u_\eps\cdot{\nabla})\theta_\eps-
\eps^{\gamma}\kappa\Delta{\theta_\eps} = {{(T_1-T_2)} \over h} (u_\eps)_n},\\[2ex]
{\rm div}\,u_\eps=0,
\end{array} \right.\;\; x\in\om,\;t\in{\bf R}^+. 
\end{equation}
For the temperature
$\theta_\eps$ we get initial data $\theta_\eps(x,0)=\theta^0_\eps(x)$
and the new boundary data
$\theta_\eps=0$ at $x_n=0$ and
at $x_n=h$. The initial and boundary data for $u_\eps$
remain unchanged. Moreover, $u_\eps$ and
$\theta_\eps$ and their gradients and ${p_\eps}$ are periodic as above.

It turns out that the value $\gamma=3/2$ is critical.
To see that we perform a multiple scales expansion technique of the unknown 
quantities
$u_\eps$, $p_\eps$ and $\theta_\eps$ and assume that
\[
u_\eps(x,t)=\eps^\rho\sum_{i=0}^{\infty}\eps^iu_i(x,{x\over\eps^\mu},t,{t\over\eps^\beta}),
\]
\[
p_\eps(x,t)=\sum_{i=0}^{\infty}\eps^ip_i(x,{x\over\eps^\mu},t,{t\over\eps^\beta}),
\]
\[
\theta_\eps(x,t)=\sum_{i=0}^{\infty}\eps^i \theta_i(x,{x\over\eps^\mu},t,{t\over\eps^\beta}),
\]
where $u_i$, $p_i$ and $\theta_i$ are all assumed to be $T^n$-periodic with respect
to $y\in{\bf R}^n$, $n=2,\,3$, $T^n$ being the usual unit torus in ${\bf R}^n$. 
If we put $y=x/\eps^\mu$ and $\tau=t/\eps^\beta$, the chain rule transforms
the differential operators as
\[
{\p\over\p{t}}\mapsto{\p\over\p{t}}+{1\over\eps^\beta}{\p\over\p{\tau}},\;\;
{\p\over\p{x}}\mapsto{\p\over\p{x}}+{1\over\eps^\mu}{\p\over\p{y}}.
\]
The question is: Can we preserve the structure from the original
problem in the leading order approximation? One
sees that if the viscosity and conductivity terms scales like $\eps^{3/2}$,
then the choices $\mu=1$ and $\rho=\beta=1/2$
preserve all quantities.
By substituting all this into the system (4.2) we can equate the
powers of $\eps$.
This formal manipulation
shows that the functions $u_i$ and $\theta_i$, $i=0,\,1,...$ are all independent of
$t$ and that $p_0$ is independent of $y$. Formally this means that
there are two scales in space and in order to preserve the evolutionary
behaviour, in the leading order approximation, a chance of time scale
$t\mapsto {t\over\eps^{1/2}}=\tau$ becomes necessary, c.f. Lions \cite{jll}.
It also follows, from the equations, that the functions $p_i$ are independent of $t$. 
The leading order system reads
\begin{equation}
\left\{ \begin{array}{l}
{\df {\p{u_0}\over\p{\tau}}+(u_0\cdot{\nabla_y})u_0
-\nu\Delta_{y}{u_0}+\nabla_y p_1 =
g \alpha e_n\theta_0}-\nabla_x p_0, \\[2ex]
{\df {\p{\theta_0}\over\p{\tau}}+(u_0\cdot{\nabla_y})\theta_0
-\kappa\Delta_{y}{\theta_0} = 0}, \\[2ex]
{\rm div_y}u_0 = 0,\;\;{\df {\rm div_x}(\int_{T^n} u_0dy)} = 0, 
\end{array} \right.
\end{equation}
on $\om\times T^n\times{\bf R}^+$
with $T^n$-periodicity in $y$.

We see that the system (4.3) differs from the system (1.4) only by the additional
forcing term coming from the "global" pressure gradient. Thus the local problem
(on the $\eps$-scale) has {\em two pressure gradients}, one local $\nabla_y p_1$ and one global
$\nabla_x p_0$.
\begin{rem}
The scaling discussed above is isotropic but the convection rolls
in Rayleigh-B\'enard convection have a preferred direction, say along the
$x_2$-axis. This commands the use of anisotropic scaling for the amplitude
equations and it turns out that the 
scaling by $\epsilon^{3/2}$ is the one pertaining to the other directions,
perpendicular to the orientation of the rolls. The details of this affine scaling will be spelled
out elsewhere.
\end{rem}

\subsection{Existence results}

In three dimensions the following global existence result holds true:
\begin{thm}
Consider the Boussinesq system (4.2). For every fixed value
of $\eps>0$ the following holds true: For every weak solution
$u_{w,\eps},\,\theta_{w,\eps}$ in
$C_w({\bf R}^+;(L^2)^4(\om))$, $\om$ a bounded open set in ${\bf R}^3$ with a large
aspect ratio, $ {L \over h} >> {g\alpha^2 (T_1-T_2)^2 L^3 \over \nu\kappa}$,
where $L$ is the length
and $h$ is the height of $\Omega$;
there exists a time $t_0$, such that there exists a unique strong
solution $u_{\eps},\,\theta_{\eps}$ in
$C([t_0,\infty);(H^1)^4(\om))$, $t\geq t_0$, with initial data
$u_{\eps}(t_0)=u_{w,\eps}(t_0)$ and $\theta_{\eps}(t_0)=\theta_{w,\eps}(t_0)$.
Moreover, the system (4.2) possesses a global attractor in $(L^2)^4(\om)$.
\end{thm}
\begin{proof}
Theorem 4.1 follows from Theorem 3.3
if we set the viscosity and heat conductivity in Theorem 3.3 equal to
$\eps^{3/2}\nu$ and $\eps^{3/2}\kappa$, respectively.
\renewcommand{\qed}{}
\end{proof}
\begin{rem}
The existence of a pressure $p_\eps\in C([t_0,\infty);(H^1)(\om)/{\bf R})$
follows by the same arguments as those in Remark 2 after the proof of Theorem 2.2.
\end{rem}

We conclude with the existence result for the two-scale system (4.3).
\begin{thm}
Consider the system (4.3) with initial data
$u_0(x,y,0)=a(x,y)$ and \linebreak[4]
$\theta_0(x,y,0)=b(x,y)$, with zero mean over $T^n$.
For almost every $x\in\om$ the system has a unique global strong solution
$u_0\in C([t_0,\infty);(H^1(T^n))^n)$, $n=2,3$, and $\theta_0\in C([t_0,\infty);H^1(T^n))$.
Moreover, $p_1\in C([t_0,\infty);H^1(T^n)/{\bf R})$. Integrated over $T^n$,
$(u_0)_i$, $\theta_0$ and $p_1$ all belong to $L^2(\om)$. Finally, 
$p_0\in C([t_0,\infty);H^1(\om))$.
In two dimensions $t_0=0$ and in three dimensions $t_0=t_0(|a|,|b|)>0$, in general.
Moreover for almost every $x\in\om$, the system (4.3) possesses a global attractor in 
$(L^2(T^n))^n\times L^2(T^n)$, $n=2,\,3$.
\end{thm}
\noindent
The proof in the two-dimensional case is relatively straightforward
and details can be found in Foias et al. \cite{fmt}.
\begin{proof}
In the three-dimensional case the proof is similar to the proof
of Theorem 3.3 but simpler. First we multiply the first equation in (4.3) by 
$\theta $ and integrate by parts to get
\[
|\theta_0(t)| \le |b|e^{-\gamma t}
\]
after integration in $t$, where $\gamma = \lambda_1 \kappa$. In \cite{bs2}
we show that
\[
g \alpha e_3\theta_0-\nabla_x p_0= g \alpha \pi_2(e_3\theta_0),
\] 
where $\pi_2(e_3\theta_0)$ denotes the projection onto the divergence free part of
$e_3\theta_0$. Hence, 
\[
|g \alpha e_3\theta_0-\nabla_x p_0|\leq g \alpha|\theta_0|\leq 
g \alpha |b|e^{-\gamma t}.
\]
This means that for $t$ sufficiently large
the forcing in the first equation of (4.3) becomes small. Thus we obtain
the existence of unique global solution
$u_0\in C([t_0,\infty);(H^1(T^n))^3)$ to the
first equation of (4.3), by Theorem 2.2, and an absorbing set in this space.
Since $u_0(t)\in (H^1(T^n))^3$ we also obtain global existence for $\theta_0$ in the
second equation of (4.3), just as in the proof of Theorem 3.3 and the existence 
of an absorbing set in this space as well.
The existence of the attractors follows from the fact
that the equations possess an absorbing set in
$(H^1(T^n))^3\times H^1(T^n)$ which is
compactly embedded in 
$(L^2(T^n))^3\times L^2(T^n)$.\\
\end{proof}
\begin{rem}
In Section 7, we prove, by using homogenization theory,
that the solutions of the scaled system (4.2) converges in the two-scale sense, 
see \cite{ng}, to the unique  solution to the  system (4.3).
\end{rem}
\begin{rem}
By averaging the system (4.3) over the unit torus in local variable
$y$ we obtain the mean-field corresponding to the scaled system (4.2) in
Section 8. The understanding
of the effects of this field on the system is crucial in the theoretical understanding
of the complex patterns in Rayleigh-B\'enard convection.
\end{rem}

\section{A priori estimates}

We are interested in the asymptotic behaviour of the
system (4.2) as $\eps\to{0}$. In order to accomplish this we need
to establish uniform (in $\eps$) bounds on $u_\eps$, $p_\eps$ and $\theta_\eps$.

This is problematic, since
the forcing term in the Navier-Stokes equation in (4.2) involves $\theta_\eps$
which can be
highly oscillatory for small values of $\eps$.
We begin with the following:
\begin{lem}
Let $u$, $p$ and $\theta$
be the solution to the system
\begin{equation}
\left\{ \begin{array}{l}
\df {{\p{u}\over\p{t}}+(u\cdot{\nabla})u-\nu\Delta{u}
+\nabla{{p}}
= e_n\theta,}\\[2ex]
\df {{\p{\theta}\over\p{t}}+(u\cdot{\nabla})\theta-
\kappa\Delta{\theta_\eps} = (u)_n},\\[2ex]
{\rm div}\,u=0,
\end{array} \right.\hfill x\in\om,\;t\in{\bf R}^+. \hfill
\end{equation}
Suppose that $u$, $p$ and $\theta$ are all
$T^n$ periodic. If
\[
\int_{T^n}u_i(y,0)dy=\int_{T^n}\theta(y,0)dy=0,
\]
then
\[
\int_{T^n}u_i(y,t)dy=\int_{T^n}\theta(y,t)dy=0,
\]
for all $t>0$.
\end{lem}
\begin{proof}
By the equation (5.1) we have
\[
\int_{T^n}({\p{u}\over\p{t}})_idy
=\int_{T^n}((e_n\theta)_i-((u\cdot{\nabla})u)_i
+\nu(\Delta{u})_i
-(\nabla{{p}})_i)dy
\]
and
\[
\int_{T^n}({\p{\theta}\over\p{t}})dy
=\int_{T^n}((u)_n-(u\cdot{\nabla})\theta
+\kappa\Delta{\theta})dy.
\]
All terms involving spatial derivatives will vanish by the
$T^n$-periodicity on the unit torus and by the incompressibility of $u$.
Moreover, $(e_n\theta)_i=0$, for $i<n$, so the system reduces to
\[
{d\over dt}\int_{T^n}(u)_n dy
=\int_{T^n}\theta dy
\]
and
\[
{d\over dt}\int_{T^n}\theta dy
=\int_{T^n}(u)_n dy.
\]
We solve this system of ODE's with the initial condition
\[
\int_{T^n}u_i(y,0)dy=\int_{T^n}\theta(y,0)dy=0,
\]
to get the
trivial solution,
\[
\int_{T^n}u_i(y,t)dy=\int_{T^n}\theta(y,t)dy=0,
\]
for all $t> 0$.
\end{proof}

We are clearly making a strong assumption assuming that the solution 
is periodic in the above Lemma. In the applications we have in mind it will
not be exactly periodic but it will be sufficiently oscillatory so that 
periodicity on a small scale is a reasonable model. Flow in porous media is
obviously one example of such a situation but others include, a mixture of
hot and cold fluid, a two fluid mixture (oil and water), a mixture of snow
and air, or water and sediment, and a turbulent fluid. This hypothesis is 
similar to the model one
uses in homogenization of materials and the test of the hypothesis is how 
well one captures averaged quantities. 

The next lemma is a Poincar\'e
inequality, c.f. L. Tartar (Lemma 1 in the Appendix of \cite{sapa}). Our version
of the lemma differs from Tartar's in that we
do not have a vanishing boundary condition on the local tori. Instead we benefit
from the result of Lemma 5.1, that the mean value vanishes.
\begin{lem}
Let $u_i\in L^2([0,s];H^1(\om))$ and assume that 
$u$ is periodic and that initially
$u$ has mean value zero over the unit torus
$T^n$, 
then
\begin{equation}
\int_0^s \int_\om|u(x,t)|^2dxdt\leq \eps^2C\int_0^s \int_\om|\nabla
u(x,t)|^2dxdt,
\end{equation}
where $C$ is a constant independent of $\eps$.
\end{lem}
\begin{proof}
We apply the Poincar\'e inequality
on $T^n$ which, by the result of Lemma 5.1, gives
\[
\int_{T^n}|u(y,t)|^2dy\leq C_1\int_{T^n}|\nabla_y u(y,t)|^2dy.
\]
A change of variables $x=\epsilon y$ yields
\[
\int_{\eps T^n}|u(x,t)|^2dx \le  \epsilon^2 C_1\int_{\eps T^n }|\nabla_x
u(x,t)|^2dx.
\]
We note that the constant $C_1$ will be the same
for all $\eps T^n$-cubes in the interior of $\om$. For $\eps T^n$-cubes
intersecting $\p\om$,
$u$ will vanish at, at least, one point and the usual Poincar\'e inequality
applies.
Let $C_2$ denote the maximum of all the constants for these cubes. 
A summation over all $\eps T^n$-cubes
gives
\[
\int_{\Omega }|u(x,t)|^2dx \le C \epsilon^2 \int_{\Omega }|\nabla_x u(x,t)|^2dx,
\]
where $C=max\{C_1,C_2\}$.
Finally, an
integration with respect to $t$ gives the desired inequality. 
\end{proof}
\begin{lem}
Let $u_\eps$ and $\theta_\eps$ satisfy (2.1) and assume that the initial data in
(4.2) is 
bounded independent of $\epsilon,$ then
\begin{equation}
\eps^{-1/2}\|u_\eps\|_{L^2([t_0,s];L^2(\om))}\leq C
\end{equation}
and
\begin{equation}
\eps^{1/2}\|\nabla u_\eps\|_{L^2([t_0,s];L^2(\om))}\leq C,
\end{equation}
for any $s>0$, where $C$ is independent of $\eps$.
\end{lem}
In the two-dimensional case
$t_0=0$ and in the three-dimensional case $t_0> 0$ in general. This
will hold true throughout the paper.
\begin{proof}
Let us consider the first equation in (4.2). We multiply by
$u_\eps$ and integrate over $\om\times]t_0,s[$. By the
incompressibility we get, by using the Schwarz's inequality
\[
{1\over 2}|u_\eps(s)|^2+\eps^{3/2}\|\nabla
u_\eps\|_{L^2([t_0,s];L^2(\om))}^2
\]
\[
\qquad\leq
\|\theta_\eps\|_{L^2([t_0,s];L^2(\om))}\|u_\eps\|_{L^2([t_0,s];L^2(\om))}+
{1\over 2}|u_\eps(t_0)|^2.
\]
By Lemma 3.1, see also the proof of Theorem 3.2,
$\|\theta_\eps\|_{L^2([t_0,s];L^2(\om))}$ is bounded, and assuming 
that $t_0$ is a time where $|\nabla u|$ is finite, see Lemma 2.1, 
we can absorb the 
term ${1\over 2}|u_\eps(t_0)|^2$ into the time integral using Lemma 5.2.
This gives the estimate,
\[
{1\over 2}|u_\eps(s)|^2+\eps^{3/2}\|\nabla
u_\eps\|_{L^2([t_0,s];L^2(\om))}^2\leq
C\|u_\eps\|_{L^2([t_0,s];L^2(\om))}
\]
so that,
\[
\eps^{3/2}\|\nabla u_\eps\|_{L^2([t_0,s];L^2(\om))}^2\leq
C\|u_\eps\|_{L^2([t_0,s];L^2(\om))}.
\]
Now we recall that $u_\eps\in L^2([t_0,s];H^1(\om))$ by Lemma 3.2
 (or Lemma 2.1). Thus, by
Lemma 5.2,
\[
\eps^{1/2}\|\nabla u_\eps\|_{L^2([t_0,s];L^2(\om))}\leq C.
\]
An application of Lemma 5.2 once again gives the desired result, i.e.,
\[
\eps^{-1/2}\|u_\eps\|_{L^2([t_0,s];L^2(\om))}\leq C.
\]
If the initial data $u_\epsilon (t_0)$ and $\theta_\epsilon (t_0)$ is bounded,
in $L^2(\Omega),$ independent of $\epsilon$
then C will also be independent of $\epsilon.$
\end{proof}

We continue with a few consequences of the above results:
\begin{cor}
Consider the first equation in (4.2). The convection term is bounded,
\begin{equation}
\|(u_\eps\cdot \nabla )u_\eps\|_{L^2([t_0,s];L^1(\om))}\leq C,
\end{equation}
where $C$ is independent of $\eps$.
\end{cor}
\begin{proof}
By the Schwarz's inequality, Lemma 3.2 and Lemma 5.3 it
immediately follows that
\[
\|u_\eps\cdot\nabla u_\eps\|_{L^2([t_0,s];L^1(\om))}\leq
\eps^{-1/2}\|u_\eps\|_{L^2([t_0,s];L^2(\om))}\eps^{1/2}\|\nabla
u_\eps\|_{L^2([t_0,s];L^2(\om))}
\le C^2.
\]
\vspace{-1.5em}
\end{proof}

\begin{cor}
Consider the first equation in (4.2). The time derivative is bounded,
\begin{equation}
\|{\p u_\eps\over \p t}\|_{L^2([t_0,s];L^2(\om))}\leq C,
\end{equation}
and the pressure is bounded,
\begin{equation}
\|p_\eps\|_{L^2([t_0,s];H^1(\om)/{\bf R})}\leq C,
\end{equation}
where $C$ is independent of $\eps$.
\end{cor}
\begin{proof}
(5.6) follows from duality by Corollary 5.1 and (5.7) follows by the
Remark 2.
\end{proof}

\section{Two-scale convergence}

In this section we recall the technically useful
concept of two-scale convergence (\cite{a2} and \cite{ng}),
for the case when the functions also depend on a time variable.
Let us consider the space $C^{\infty}(T^n)$ of smooth periodic functions, with
unit period, in ${\bf R}^n$.
\begin{definition}
A sequence $\{u_\eps\}$ in $L^2([0,s];L^2(\om))$ is said to two-scale converge
to $u_0=u_0(x,y,\tau)$ in $L^2([0,s];L^2(\om\times T^n))$ if, for
any
$\varphi\in C^\infty_0(\om\times[0,s];C^{\infty}(T^n))$,
\begin{equation}
\int_0^s\int_\om u_\eps(x,\tau)\varphi(x,{x\over \eps},\tau)dxd\tau\to
\int_0^s\int_\om\int_{T^n} u_0(x,y,\tau)\varphi(x,y,\tau)dydxd\tau,
\end{equation}
as $\eps\to 0$.
\end{definition}
\noindent
We have the following extension of a compactness result first proved by
Nguetseng \cite{ng} and then further developed by Allaire in \cite{a1,a3,a2}
and by Holmbom in \cite{ho}.
\begin{thm}
Suppose that $\{u_\eps\}$ is a uniformly bounded sequence in
$L^2([0,s];L^2(\om))$.
Then there exists a subsequence, still denoted by $\{u_\eps\}$, and a function
$u_0=u_0(x,y,\tau)$ in $L^2([0,s];L^2(\om\times T^n))$, such that
$\{u_\eps\}$ two-scale converges to $u_0$.
\end{thm}
\noindent
The relation between ${\tilde u}$ and $u_0$ is explained in the next theorem.
In fact, by choosing test functions which do not depend
on $y$ we  have the following (see e.g. \cite{a1}):
\begin{thm}
Suppose that $\{u_\eps\}$ two-scale converges to $u_0$, where
$u_0=u_0(x,y,\tau)$ in $L^2([0,s];L^2(\om\times T^n))$, then
$\{u_\eps\}$ converges to ${\overline {u_0}}$ weakly in $L^2([0,s];L^2(\om))$,
where
\[
{\overline {u_0}}(x,\tau) = \int_{T^n} u_0(x,y,\tau)dy.
\]
\end{thm}
In other words:
\[
{\tilde u}={\overline {u_0}}.
\]
\begin{rem}
The results of Theorem 6.1 and Theorem 6.2 remain valid for the
larger class of {\em admissible} test functions
$\varphi\in L^2([0,s]\times\om;C(T^n))$, see e.g. \cite{a1}. 
\end{rem}

\section{Homogenization}

With the help of the results from Section 4, Section 5 and the 
two-scale compactness result Theorem 6.1 from the previous section
we can now state the main results of the paper. 
\begin{thm} Consider the Navier-Stokes system (4.2). Suppose that
the initial data
$u_\eps(x,0)=u^0_\eps(x)$ and $\theta_\eps(x,0)=\theta^0_\eps(x)$
two-scale converge to unique limits
$u^0(x,y)$ and $\theta^0(x,y)$, respectively.
Then, as $\eps\to 0$, the following quantities two-scale converge,
\[
\eps^{-1/2}u_\eps\to u_0,
\]
\[
p_\eps\to p_0,
\]
\[
\theta_\eps\to \theta_0,
\]
where $p_0=p_0(x,\tau)$.
Moreover, there exists a function $p_1=p_1(x,y,\tau)$
such that
\[
\nabla p_\eps\to \nabla_x p_0+\nabla_y p_1,
\]
the functions 
$u_0$, $p_0$, $p_1$ and $\theta_0$,
being the unique solutions to the Navier-Stokes system (4.3)
with  initial data
$u_0(x,y,0)=u^0(x,y)$, $\theta_0(x,y,0)=\theta^0(x,y)$ and boundary data
$T^n$-periodic in the variable $y$.
\end{thm}
\noindent
Before we prove the theorem we state a corollary
\begin{cor}
Consider the Navier-Stokes system (4.2).
The following quantities two-scale converge,
\[
\eps^{-1/2}u_\eps\to u_0,
\]
\[
p_\eps\to p_0,
\]
\[
\nabla p_\eps\to \nabla_x p_0+\nabla_y p_1,
\]
\[
T_\eps\to T_0,
\]
where
$u_0$, $p_0$, $p_1$ and $T_0$ are the unique solution to
the Navier-Stokes system (4.3).
\end{cor}

\begin{proof}[Proof of Theorem 7.1]
By assumption the initial data admit unique two-scale limits.
Consider now the first equation of (4.2),
\[
\df
{{\p{u_\eps}\over\p{t}}+(u_\eps\cdot{\nabla})u_\eps-\eps^{3/2}\nu\Delta{u_\eps}
+\nabla{{p_\eps}}
= e_2\theta_\eps,\;\mbox{ in }\om\times{\bf R}^+}.
\]
Let $s>0$ and choose test functions
$\varphi\in C^\infty_0(\om\times[t_0,s];C^{\infty}(T^n))$.
By the results
of Section 5, all the sequences $\{\eps^{-1/2}u_\eps\}$,
$\{{\p u_\eps\over \p t}\}$, $\{p_\eps\}$
and $\{\theta_\eps\}$ are uniformly bounded in $L^2([t_0,s];L^2(\om))$.
Therefore, according to Theorem 6.1, they admit two-scales limits.
In order to identify these limits we multiply
each of these terms by the smooth compactly supported test function
$\varphi(x,{x\over \eps},\tau)$. For the time derivative we get (recall that
$\tau=t/\sqrt{\eps}$)
\[
\int_0^s\int_\om {\p{u_\eps}\over\p{t}}
(x,t)\varphi(x,{x\over\eps},\tau)dxd\tau=
\eps^{-1/2}\int_0^s\int_\om u_\eps{\p\varphi\over\p\tau}dxd\tau.
\]
Sending $\eps\to 0$ and integrating by parts yield,
\[
\int_0^s\int_\om\int_{T^n}{\p{u_0}\over\p{\tau}}\varphi(x,y,\tau)dydxd\tau.
\]
For the second (inertial) term we consider
\[
\int_0^s\int_\om [(u_\eps\cdot{\nabla}u_\eps(x,\tau)-u_0\nabla_y 
u_0(x,y,\tau)]
\varphi(x,{x\over \eps},\tau)dxd\tau
\]
\[
\qquad=\int_0^s\int_\om [\eps^{-1/2}u_\eps\cdot(\eps^{1/2}\nabla u_\eps-\nabla_y u_0)
+(\eps^{-1/2}u_\eps-u_0)\nabla_y u_0]
\varphi(x,{x\over \eps},\tau)dxd\tau.
\]
By considering $\nabla_y u_0\varphi(x,{x\over \eps},\tau)$ to be a
test function in the second term
on the right hand side
this term immediately passes to zero in the two-scale sense. For the first term
the
Schwarz's inequality yields
\[
\int_0^s\int_\om\eps^{-1/2}u_\eps\cdot(\eps^{1/2}\nabla u_\eps-\nabla_y 
u_0)\cdot
\varphi(x,{x\over \eps},\tau)dxd\tau
\]
\[
\qquad\leq
\|\eps^{-1/2}u_\eps\|\|(\eps^{1/2}\nabla u_\eps-\nabla_y u_0)
\varphi(x,{x\over \eps},\tau)\|
\]
\[
\qquad\leq{C}\|(\eps^{1/2}\nabla u_\eps-\nabla_y u_0)
\varphi(x,{x\over \eps},\tau)\|,
\]
where all the norms are in ${L^2([0,s];L^2(\om))}$.
In order to pass to the limit in the right hand side we consider the usual
$L^2$-mollifications $\nabla u_{\eps,\mu}$ and $\nabla_y u_{0,\mu}$ of
$\eps^{1/2}\nabla u_\eps$ and $\nabla_y u_{0}$, respectively.
Since the mollified functions pass strongly to
$\eps^{1/2}\nabla u_\eps$ and $\nabla_y u_{0}$, respectively,
as $\mu\to 0$, we have, for $\mu$ sufficiently small, say $\mu\leq \mu_0$ and
for every $y$ in $T^n$   
\[
\|(\eps^{1/2}\nabla u_\eps-\nabla_y u_0)\varphi\|
\leq\|(\nabla u_{\eps,\mu}-\nabla_y u_{0,\mu})\varphi\|+\delta,
\]
where $\delta$ is arbitrarily small, independently of $\eps$. This inequality
still holds
true if we take the supremum in $y$ over $T^n$.
Thus, for every $\mu\leq \mu_0$, the right hand side
will tend to zero as $\eps$ tends to $0$, by the uniqueness of the
two-scale limit $\nabla_y u_0$ of $\eps^{1/2}\nabla u_\eps$.
Consequently, we have proved that 
\[
u_\eps\nabla u_\eps\to u_0\nabla_y u_0
\]
in the two-scale sense.
For the third term we get, by the divergence theorem,
\[
-\eps^{3/2}\nu\int_0^s\int_\om\Delta{u_\eps}(x,\tau)\varphi(x,{x\over
\eps},\tau)dxd\tau
\]
\[
\qquad=-\eps^{-1/2}\nu\int_0^s\int_\om u_\eps(x,\tau)\Delta_{y}\varphi(x,{x\over
\eps},\tau)dxd\tau+
\;\;{\rm terms\; tending\; to\; zero\; as}\;\;\eps\to 0.
\]
Sending $\eps\to 0$ yields, after applying the divergence theorem again,
\[
-\nu\int_0^S\int_\om\int_{T^n}\Delta_{y}{u_0}\varphi(x,y,\tau)dydxd\tau.
\]
For the right hand side we immediately get
\[
\int_0^s\int_\om e_2\theta_\eps(x,\tau)\varphi(x,{x\over \eps},\tau)dxd\tau\to
\int_0^s\int_\om\int_{T^n} e_2\theta_0\varphi(x,y,\tau)dydxd\tau. 
\]
For the fourth term (the pressure) we have to be a bit more careful.
Let us multiply the first equation of (4.2) by $\eps\varphi(x,{x\over
\eps},\tau)$.
For the pressure term we get
\[
\eps\int_0^s\int_\om\nabla{p_\eps}(x,\tau)\varphi(x,{x\over \eps},\tau)dxd\tau
=\int_0^s\int_\om{p_\eps}(x,\tau)(\eps{\rm div}_x+{\rm div}_y)\varphi(x,{x\over
\eps},\tau)
dxd\tau.
\]
A passage to the limit, and an application of the divergence theorem, using the
fact that
all other terms vanish, gives
\[
\int_0^s\int_\om\int_{T^n}\nabla_y{p_0}(x,y,\tau)\varphi(x,y,\tau)dydxd\tau=0,
\]
which implies that $p_0$ does not depend on $y$.
We now add the local incompressibility assumption on the test functions
$\varphi$,
i.e. ${\rm div}_y\varphi=0$ and multiply the pressure term by $\varphi(x,{x\over
\eps},\tau)$
and apply the divergence theorem,
\[
\int_0^s\int_\om\nabla{p_\eps}(x,\tau)\varphi(x,{x\over \eps},\tau)dxd\tau
=\int_0^s\int_\om{p_\eps}(x,\tau){\rm div}_x\varphi(x,{x\over \eps},\tau)
dxd\tau.
\]
A passage to the limit, and an application of the divergence theorem, gives
\[
\int_0^s\int_\om\int_{T^n}\nabla_x{p_0}(x,\tau)\varphi(x,y,\tau)dydxd\tau.
\]

Collecting all two-scale limits on the right hand side gives
\[
\int_0^s\int_\om\int_{T^n}(f-{\p{u_0}\over\p{\tau}}-(u_0\cdot{\nabla}_y)u_0+\nu\Delta_{y}{u_0}
-\nabla_x{p_0})\varphi dydxd\tau = 0.
\]
Since ${\rm div}_y\varphi=0$
we can argue as in Remark 2 and conclude that 
there exists
a local pressure gradient $\nabla_y p_1(x,y,\tau)$ given by
\[
\nabla_y p_1(x,y,\tau)
=f-{\p{u_0}\over\p{\tau}}-(u_0\cdot{\nabla}_y)u_0+\nu\Delta_{y}{u_0}
-\nabla_x{p_0}.
\]
Let us now consider the second equation of (4.2). We already know that the
sequence
$\{\theta_\eps\}$ is uniformly bounded in $L^2([0,s];L^2(\om))$.
We multiply by $\eps^{1/2}\varphi$ as above and
for the time derivative we get
\[
\eps^{1/2}\int_0^s\int_\om{\p{\theta_\eps}\over\p{t}}\varphi dxd\tau=
\int_0^s\int_\om \theta_\eps{\p\varphi\over\p\tau}dxd\tau.
\]
By letting $\eps\to\ 0$ we get
\[
\int_0^s\int_\om \theta_\eps{\p\varphi\over\p\tau}dxd\tau\to
\int_0^s\int_\om\int_{T^n} {\p\theta_0\over\p\tau}\varphi dydxd\tau.
\]
For the non-linear term we have
\[
\eps^{1/2}\int_0^s\int_\om(u_\eps\cdot{\nabla})\theta_\eps\varphi dxd\tau
\]
\[
\qquad=\eps^{1/2}\int_0^s\int_\om(u_\eps\cdot{\nabla}_x)\varphi\theta_\eps dxd\tau+
\eps^{-1/2}\int_0^s\int_\om(u_\eps\cdot{\nabla}_y)\varphi\theta_\eps dxd\tau.
\]
The first term on the right hand side immediately passes to zero. For the second
term we
argue as above and consider the difference
\[
\int_0^s\int_\om((\eps^{-1/2}u_\eps\cdot{\nabla}_y)\varphi\theta_\eps-
(u_0\cdot{\nabla}_y)\varphi\theta_0) dxd\tau
\]
\[
\qquad=\int_0^s\int_\om((\eps^{-1/2}u_\eps\cdot{\nabla}_y)\varphi\theta_\eps-
(u_0\cdot{\nabla}_y)\varphi\theta_\eps) dxd\tau
\]
\[
\qquad\quad+\int_0^s\int_\om((u_0\cdot{\nabla}_y)\varphi\theta_\eps-
(u_0\cdot{\nabla}_y)\varphi\theta_0) dxd\tau.
\]
By considering $(u_0\cdot{\nabla}_y)\varphi$ to be a test function in the second
term,
this term immediately passes to zero in the two-scale sense. For the first term
we get, by the Schwarz's inequality,
\[
\int_0^s\int_\om((\eps^{-1/2}u_\eps\cdot{\nabla}_y)\varphi\theta_\eps-
(u_0\cdot{\nabla}_y)\varphi\theta_\eps) dxd\tau\leq
{C}\|(\eps^{-1/2}u_\eps-u_0)\cdot{\nabla}_y\varphi\|,
\]
where the norm is in ${L^2([0,s];L^2(\om))}$.
We introduce, as above, mollifiers and consider the sequence
$\{u_{\eps,\mu}\}$ which converges to $\eps^{-1/2}u_\eps$ strongly as
$\mu\to 0$. Arguing as above, we choose $\mu$ sufficiently small
to get
\[
\|(\eps^{-1/2}u_\eps-u_0)\cdot{\nabla}_y\varphi\|
\leq \|(u_{\eps,\mu}-u_0)\cdot{\nabla}_y\varphi\|+\delta,
\]
where $\delta$ is arbitrarily small. Thus, by sending $\eps\to 0$,
\[
u_\eps\nabla \theta_\eps\to u_0\nabla_y \theta_0
\]
in the two-scale sense.
For the third term we get, by the divergence theorem,
\[
-\eps^{2}\kappa\int_0^s\int_\om\Delta{\theta_\eps}\varphi dxd\tau
\]
\[
\qquad=-\kappa\int_0^s\int_\om \theta_\eps\Delta_{y}\varphi dxd\tau+
\;\;{\rm terms\; tending\; to\; zero\; as}\;\;\eps\to 0.
\]
We let $\eps\to 0$ and get, by the divergence theorem,
\[
-\kappa\int_0^s\int_\om \theta_\eps\Delta_{y}\varphi dxd\tau\to
-\kappa\int_0^s\int_\om\int_{T^n} \Delta_{y}\theta_0\varphi dydxd\tau.
\]
The right hand side of the second equation in (4.2) will vanish 
since $u_\eps$ is of order $\eps^{1/2}$.
By  Theorem 4.2
the system
(4.3) has a unique solution $\{u_0,p_0,p_1,\theta_0\}$ and, thus, by
uniqueness,
(4.3) is two scales
homogenized limit of the system (2.1).
Also, by uniqueness, the whole sequence converges
to its two-scale limit and the theorem is proven.
\end{proof}

\section{The mean velocity field}

In this section we derive the mean field ${\overline {u}}_0$ for the velocity.
Let us consider the first equation of (4.3)
\[
\left\{ \begin{array}{l}
{\df {\p{u_0}\over\p{\tau}}+(u_0\cdot{\nabla_y})u_0
-\nu\Delta_{y}{u_0}+\nabla_y p_1 =
e_n\theta_0}-\nabla_xp_0, \\[2ex]
{\df {\rm div_y}u_0 = 0,\;\;{\rm div_x}(\int_{T^n}\,u_0\,dy) = 0}, 
\end{array} \right.
\]
in $\om\times {T^n}\times{\bf R}^+$
with $T^n$-periodicity as boundary data in $y$.

By letting $K=K(y,y',\tau)$ denote
the heat kernel we can write the solution to
the first equation of (4.3) as an integral,
\[
u_0(x,y,\tau)=\int_{T^n}K(y,y',\tau)u_0(x,y',0)dy'
\]
\[
\qquad+\int_0^\tau\int_{T^n}K(y,y',\tau-s)(e_n\theta_0(x,y',s)-\nabla_xp_0(x,s)-\nabla_y
p_1(x,y',s)
\]
\[
\qquad\qquad\qquad-(u_0\cdot{\nabla_y})u_0(x,y',s))dy'ds.
\]

Now, by the decay of the heat kernel, for $\tau$ sufficiently large,
the first term becomes arbitrarily small. An averaging of the second term
over $T^n$ in $y$ gives
\[
{\overline {u_0}}(x,\tau)=
\int_0^\tau(e_n{\overline {\theta_0}}-\nabla_xp_0-
{\overline {(u_0\cdot{\nabla_y})u_0}})(x,s)ds.
\]
The divergence theorem gives
\[
\int_{T^n}(u_0\cdot \nabla_y)u_0\,dy=-\int_{T^n}({\rm div}_y u_0)u_0dy
+\int_{\p T^n}(u_0\cdot n_y)u_0\,dS_y,
\]
where $n_y$ is the local unit normal and $S_y$ the local surface element.
Now, by the local incompressibility we get
\[
\int_{T^n}({\rm div}_y u_0)u_0dy=0.
\]
Moreover, for the boundary integral we get, in the two-dimensional case,
\[
\int_{\p T^2}u_in_iu_0\,dS_y
=\int_0^1-u_2u_0(y_1,0)dy_1+u_1u_0(1,y_2)dy_2
\]
\[
\qquad+\int_0^1u_2u_0(y_1,1)dy_1-u_1u_0(0,y_2)dy_2
=0,
\]
by the $T^2$-periodicity of $u_0=(u_1,u_2)$. Consequently
\[
{\overline {(u_0\cdot{\nabla_y})u_0}}=0.
\]
The computation in the three-dimensional case is similar.
Therefore the mean field reduces to
\begin{equation}
{\overline {u_0}}(x,\tau)=
\int_0^\tau(e_n{\overline {\theta_0}}-\nabla_xp_0)(x,s)ds.
\end{equation}
We can now use the incompressibility, which when applied to (8.1) gives
\begin{equation}
{\rm div}_x(e_n{\overline \theta_0}(x,\tau)-
\nabla_xp_0(x,\tau))=0.
\end{equation}
This says that
\begin{equation}
e_n{\overline \theta_0}(x,\tau)-
\nabla_xp_0(x,\tau)=H(x,\tau)
\end{equation}
where $H$ is a divergence free (rotational) field.

The field $H$ can be determined explicitly
by solving the global equation (8.2), (in $x$), for the pressure $p_0$.
For this purpose one needs to impose the appropriate boundary condition
on the pressure in order to close the system. We impose the Neumann condition
\[
n\cdot{\nabla_x}p_0=0
\]
for the pressure. In fact
Lions \cite{jll} and Sanchez-Palencia
\cite{sapa} impose the condition
\[
n\cdot{\overline {u_0}}=0
\]
which when inserted in (8.1) gives
\[
n\cdot{\nabla_x}p_0=0.
\]

We consider
\begin{equation}
\left\{ \begin{array}{l}
\Delta_{x}p_0(\cdot,\tau)={\rm div}_x\,
{e_n{\overline \theta_0}}(\cdot,\tau),\;\;{\rm in}\;\;\om,\\[2ex]
n\cdot{\nabla_x}p_0(\cdot,\tau)=0,\;\;{\rm on}\;\;\p\om, 
\end{array} \right.
\end{equation}
for every $\tau\geq 0$,
and obtain
\[
H=\pi (e_n{\overline \theta}_0)=-\nabla\times(\Delta^{-1}(\nabla\times
e_n{\overline \theta}_0)),
\]
where $\pi (e_n{\overline \theta}_0)$ denotes projection onto the divergence
free part
of the conduction term.

Collecting the results
from the above discussion we can now express the mean field as
\[
{\overline u}_0(x,{t\over\sqrt{\eps}})=\int_0^{t/\sqrt{\eps}}
\pi (e_n{\overline \theta}_0)(x,s)ds.
\]
This gives the contribution of the
conduction and the global pressure to the small scale flow.
We have let the local flow settle down and averaged (in $y$), denoted by
overbar, over the unit cell $T^n$, $n=2,\,3$.
\begin{rem}
If we take boundary layer effects into account, we can,
considering
the simplest case of a boundary layer, specify the value of the global pressure
gradient
$\nabla_x p_0$ at the boundary, of the boundary layer. For instance
we can put
\[
n\cdot{\nabla_x}p_0=f.
\]
This will result in the additional term $G*f$ in the mean field, where
$G$ is the usual Neumann kernel for the Laplacian. In this case the
mean field becomes
\[
{\overline u}_0(x,{t\over\sqrt{\eps}})=\int_0^{t/\sqrt{\eps}}
\left(\pi (e_n{\overline \theta}_0)+G*f\right)(x,s)ds.
\]
\end{rem}
\begin{rem}
The above formulas give the mean field flow with the global
convection (rolls)
taken out. This was done by imposing the periodic boundary conditions
on the local cell, see the introduction.
The two-scale convergence carries over to a larger class of test functions which are non-periodic
, see \cite{ho}.
This will not be repeated here, but to be able to compare
the true (collective) mean field we compute the average
over the convection term in the non-periodic case. By the Taylor and mean value
theorems
we get, in the two-dimensional case,
\[
\int_{\p T^2}u_in_iu_0\,dS_y
=\int_0^1(-u_2u_0(y_1+Y_1,Y_2)+u_2u_0(y_1+Y_1,1+Y_2))dy_1
\]
\[
\qquad\quad+\int_0^1(u_1u_0(1+Y_1,y_2+Y_2)-u_1u_0(Y_1,y_2+Y_2))dy_2
\]
\[
\qquad={\overline u_2}\int_0^1(u_0(y_1+Y_1,Y_2)-u_0(y_1+Y_1,1+Y_2))dy_1
\]
\[
\qquad\quad+{\overline u_1}\int_0^1(u_0(1+Y_1,y_2+Y_2)-u_0(Y_1,y_2+Y_2))dy_2
+O(\Delta)
\]
\[
\qquad={\overline u_2}{\p{\overline u_0}\over\p Y_2}+{\overline u_1}{\p{\overline
u_0}\over\p Y_1}+
O(\Delta)={\overline u_0}\cdot{\nabla \overline {u_0}}+
O(\Delta),
\]
where $\Delta=(y_1-{\overline y_1},y_2-{\overline y_2})$,
$({\overline y_1},{\overline y_2})$ is the point in $T^2$ where u attains its
mean value, and $(Y_1,Y_2)$ is the location of the box. 
The error term $O(\Delta)\leq C\|u_0\||\Delta|$ is bounded by Theorem 4.2.
In global coordinates we therefore get
\[
{\overline {u_0\cdot\nabla u_0}}={\overline u_0}\cdot{\nabla \overline {u_0}}+
O(\eps),
\]
and similarly
\[
{\overline {\nabla p}_1} = \nabla {\overline p}_1.
\]
If we insert this into the expression for the mean field and differentiate with
respect to the fast time variable $\tau$, we obtain, as $\eps\to 0$,
\[
{\p{\overline u_0}\over \p \tau}+{\overline u_0}\cdot{\nabla \overline
{u_0}}+\nabla {\overline p}_1=
\pi (e_2{\overline \theta}_0),
\]
i.e. a forced Euler equation, where $u_0$ is a function of the scaled
variables $(\tau, Y) = (t/\sqrt{\eps}, x/\eps)$.
The computation in the three-dimensional case is similar.
Thus the local flow satisfies the Navier-Stokes equation, whereas the mean 
flow satisfies the Euler equation.
\end{rem}

\subsection*{Acknowledgements}

The first author was partially supported by NSF grants
DMS-9704874, BCS-9819095, a UCSB faculty research grant and a grant from Rann\'{\i}s and the second author was partially 
supported by a postdoctoral fellowship from the Swedish Research Council for 
the Natural Sciences.

\footnotesize

\label{bir_svan_lp}
 
\end{document}